%% file: jumps.tex
\documentclass[12pt]{article}    
\usepackage{array}
\usepackage{url}
\usepackage{amsmath,amsfonts,amssymb}
\newenvironment{Proof}{
{\setlength{\parindent}{0.0in}

{\em Proof.} }
}{
\rule[-2.2mm]{1.5mm}{3.7mm}  
\vspace{0.13in} }
\newenvironment{Sketch}{
{\setlength{\parindent}{0.0in}

{\em Sketch of proof.} }
}{
\rule[-2.2mm]{1.5mm}{3.7mm}  
\vspace{0.13in} }
\def\veestrut{\rule{0mm}{3mm}}
\def\mm{\,\o{\veestrut\wedge}\,}\def\jj{\,\o{\veestrut\vee}\,}%
\newtheorem{theorem}{Theorem}
\newtheorem{lemma}[theorem]{Lemma}
\newtheorem{corollary}[theorem]{Corollary}

\newcommand{\Dash}{\rule[.9mm]{1.5cm}{.1mm}\hspace{2mm}}
\newcommand{\EQ}{\mbox{$\;\; = \;\;$}}
\newcommand{\Wavy}{\;\approx\;}
\newcommand{\WAVY}{\;\approx\;}

\newcommand{\Compos}{\!\circ\! }
\newcommand{\FROM}{\!:\!}
\newcommand{\TO}{\longrightarrow}
\newcommand{\GOESTO}{\longmapsto}

\newcommand{\la}{\langle}
\newcommand{\ra}{\rangle}
\newcommand{\Reals}{{\mathbb R}}

\newcommand{\Integers}{{\mathbb Z}}

\newcommand{\ENDPROOF}{\hspace{0.03in}
\rule[-2.2mm]{1.5mm}{3.7mm}  
\vspace{0.03in}\\}\def\o{\overline}

\def\wavy{\approx}
\def\Models{\;\models\;}

\def\o{\overline}

\newcommand{\Tri}[3]{\mbox{$\bigtriangleup#1#2#3$}} 
\def\join{\vee}     
\def\meet{\wedge}   

\begin{document}
\begin{center}
    {\bf Discontinuities in the identical satisfaction
       of equations}\\[0.1in]
                   by\\[0.1in]
      {\bf Walter {
      {Taylor}}} (Boulder, CO)
\end{center}

\begin{center}{\bf Preliminary report, \today.}\\[0.1in]
\end{center} \vspace{0.1in}

{\small {\bf Abstract.} For a metric space $(A,d)$, and a set
$\Sigma$ of equations, some quantities are introduced that measure
the size of discontinuities that must occur in operations
satisfying $\Sigma$ (identically) on $A$.  We are able to evaluate
these quantities in a few easy cases.} \vspace{0.05in}

\tableofcontents \newpage

 \setcounter{section}{-1}

\section{Introduction.}                  \label{sec-intro}

This paper is part of a continuing investigation---see the
author's papers \cite{wtaylor-cots} (1986), \cite{wtaylor-sae}
(2000), and \cite{wtaylor-eri} (2006)---into the {\em
compatibility} relation (see (\ref{eq-space-models-sigma}) below)
between a topological space $A$ and a set $\Sigma$ of equations,
which we will briefly review in \S\ref{sub-compat-backg}.

The main results so far are Theorems  \ref{thm-23},
\ref{th-idem-comm-23}, \ref{th-s1-majority},
\ref{thm-mu-s2-majority} and \ref{th-triode}, in
\S\ref{sub-mult-zero-one}, \S\ref{sub-idem-comm},
\S\ref{sub-majority}, \S\ref{sub-s2} and \S\ref{sub-Y}. They
showcase a new quantity that measures the
in\-com\-pati\-bility---see \S\ref{sub-compat-backg}---of a metric
space $A$ with an equational theory $\Sigma$. This quantity will
be denoted $\mu(A,\Sigma)$---see \S\ref{sub-cont-fail}. In the
first four of the aforementioned theorems we are able to calculate
some non-trivial $\mu$-values, with $A=S^1$, the one-dimensional
sphere, and $\Sigma$ taken, for example, as the well-known ternary
majority laws. Then Theorem \ref{th-triode} deals with lattice
theory on a $Y$-shaped space. It invokes some more sophisticated
measures, $\mu_n$ and $\mu_n^{\star}$, which measures jumps in
$n$-fold iterates of the operations.

These calculations indicate the feasibility of studying this
$\mu$, and show that the quantity may have some independent
interest. At the same time, we seek a broader applicability of the
concept and associated methods. Thus we will for a while keep the
project open, both for the ultimate form of this paper, and for
future writings.

At the start of \S\ref{sub-cont-fail} we contrast this $\mu$ with
an earlier measure \cite{wtaylor-asi} of incompatibility known as
$\lambda$.

\subsection{Compatibility---context and background.}
            \label{sub-compat-backg}
In this context, $\Sigma$ typically denotes a set (finite or
infinite) of equations\footnote{%
A (formal) equation is an ordered pair of terms $(\sigma,\tau)$,
more frequently written $\sigma\wavy\tau$. As such it makes no
assertion, but merely presents two terms for consideration. The
actual mathematical assertion is made by the satisfaction relation
$\models$.}, %
which are understood as universally quantified. We usually expect
that $\Sigma$ has a specified similarity type. This means that we
are given a set $T$ and whole numbers $n_t\geq0$ ($t\in T$), that
for each $t\in T$ there is an operation symbol\footnote{%
In examples we may sometimes give the operation symbols familiar
names like $+$ or $\wedge$, or use $F$, $G$, etc.\ without a
subscript. All of these variations may be thought of as colloquial
expressions
for the more formal $F_t$.} %
$F_t$ of arity $n(t)$, and that the operation symbols of $\Sigma$
are included among these $F_t$.

Given a set $A$ and for each $t\in T$ a function $\o{F_t}\FROM
A^{n(t)} \TO A$ (called an operation), we say that the operations
$\o{F_t}$ {\em satisfy} $\Sigma$ and write
\begin{equation}             \label{eq-opening-set}
        (A,\o{F}_t)_{t\in T} \Models \Sigma,
\end{equation}
if for each equation $\sigma\wavy\tau$ in $\Sigma$, both $\sigma$
and $\tau$ evaluate to the same function when the operations
$\o{F_t}$ are substituted for the symbols $F_t$ appearing in
$\sigma$ and $\tau$. Given a {\em topological space} $A$ and a set
of equations $\Sigma$, we write
\begin{equation}             \label{eq-space-models-sigma}
        A \Models \Sigma,
\end{equation}
and say that $A$ and $\Sigma$ are {\em compatible}, iff there
exist {\em continuous} operations $\o{F_t}$ on $A$ satisfying
$\Sigma$.

While the definitions are simple, the relation
(\ref{eq-space-models-sigma}) remains mysterious. The algebraic
topologists long knew that the $n$-dimensional sphere $S^n$ is
compatible with H-space theory ($x\cdot e\Wavy x\Wavy e\cdot x$)
if and only if $n=1,3$ or $7$. For $A=\Reals$, the relation
(\ref{eq-space-models-sigma}) is algorithmically undecidable for
$\Sigma$ \cite{wtaylor-eri}; i.e.\ there is no algorithm that
inputs an arbitrary finite $\Sigma$ and outputs the truth value of
(\ref{eq-space-models-sigma}) for $A=\Reals$. In any case,
(\ref{eq-space-models-sigma}) appears to hold only sporadically,
and with no readily discernable pattern.

The mathematical literature contains many scattered examples of
the truth or falsity of specific instances of
(\ref{eq-space-models-sigma}). The author's earlier papers
\cite{wtaylor-cots}, \cite{wtaylor-sae}, \cite{wtaylor-eri},
\cite{wtaylor-asi} collectively refer to most of what is known,
and in fact many of the earlier examples are recapitulated
throughout the long article \cite{wtaylor-asi}. We will therefore
not attempt to write a list of examples for this introduction.

\subsection{Metric approximation to compatibility.}
              \label{sub-metric-approx-compat}

If the topological space $A$ is metrizable, and if a metric $d$ is
selected for $A$, then in addition to the modeling relation
$\models$, one may also study some real-valued measures of
approximate satisfaction.  In our previous paper
\cite{wtaylor-asi} we wrote
\begin{align}                 \label{eq-models-varep}
              A \models_{\varepsilon} \Sigma
\end{align}
(for real $\varepsilon>0$) to mean that there exist continuous
operations $\o F_t$ on $A$ such that, for each equation
$\sigma\wavy\tau$ in $\Sigma$, the terms $\sigma$ and $\tau$
evaluate to functions $\o \sigma$ and $\o \tau$ that are within
$\varepsilon$ of each other. We then defined $\lambda_A(\Sigma)$
to be {\em the smallest non-negative\/\footnote{%
If there is any real number satisfying this condition, then there
is a smallest one, by completeness. If there is no such real
number, then we let $\lambda_A(\Sigma)=\infty$.} %
real such that $A \models_{\varepsilon} \Sigma$ for every
$\varepsilon>\lambda_A(\Sigma)$.}

As one might imagine, the precise value of $\lambda_A(\Sigma)$
depends strongly on the metric $d$ chosen to represent the
topology of $A$; moreover, its value can increase if $\Sigma$ is
augmented by the inclusion of some of its own logical
consequences. The earlier article \cite{wtaylor-asi} illustrates
these points with detailed estimations of $\lambda_A(\Sigma)$ for
many different $A$ and $\Sigma$.

\subsection{Measuring continuity-failure in models of $\Sigma$.}
              \label{sub-cont-fail}

In looking to have $A\models \Sigma$ for a space $A$ and a set
$\Sigma$ of equations, we  demand both the continuity of the
operations $\o F_t$ of $(A,\o{F}_t)_{t\in T}$, and the exact
satisfaction of $\Sigma$ by these operation. The outlook reviewed
in \S\ref{sub-metric-approx-compat} was to relax the need for
exact satisfaction, and to see how close we can come with
approximate satisfaction.

In this paper we examine a different---opposite, really---way of
relaxing our requirements. Namely, we require exact satisfaction
while measuring how far our operations must deviate from
continuity.


%

Let $B$ be a topological space, and $(A,d)$ a metric space. Let us
consider a function
 $\o F\FROM B\TO A$. If $\o F$ is
continuous at $b\in B$, then  for each $\varepsilon>0$ there is a
neighborhood $U$ of $b$ with $\o F[U]\subseteq B(\o
F(b),\varepsilon/2)$ (the open ball about $\o F(b)$ with radius
$\varepsilon/2$). Consequently, if $\o F$ is continuous at $b\in
B$, then for each $\varepsilon>0$ there is a neighborhood $U$ of
$b$ such that $\o F[U]$ has diameter $<\varepsilon$, in other
words
\begin{align}                    \label{eq-def-chi-0}
            \inf \,\{\, \text{diameter}\,\o F[U]\,:\,
                              U \text{ open},\; b\in U \,\} \EQ 0.
\end{align}

If $\o F$ is not continuous at $b$ we may still define the
quantity
\begin{align}                     \label{eq-def-chi}
           \chi(\o F,b) \EQ  \inf \,\{\, \text{diameter}\,\o F[U]\,:\,
                              U \text{ open},\; b\in U \,\},
\end{align}
which we will call the {\em jump of $\o F$ at $b$}, as a real
number (or $\infty$), and take this quantity as a measure of
failure of continuity at $b$. The following lemma is almost
immediate.

\begin{lemma}          \label{lem-continuity}
$\o F\FROM B\TO A$ is continuous at $b\in B\,$ iff $\,\chi(\o
F,b)=0$.
\end{lemma}\begin{Proof}
We already saw that if $\o F$ is continuous at $b$, then $\chi(\o
F,b)=0$. For the converse, let us be given  $\chi(\o F,b)=0$ and
prove continuity at $b$. Suppose we are given $\varepsilon>0$. The
infimum appearing in (\ref{eq-def-chi-0}) and (\ref{eq-def-chi})
is zero, which means in particular that there is an open
$U\subseteq B$ with $b\in U$ and such that $\o F[U]$ has diameter
$<\varepsilon$. In other words, $d(f(u),f(b))<\varepsilon$ for all
$u\in U$. Thus $\o F$ is continuous at $b$.
\end{Proof}

Then $\chi(\o F)$, the jump of $\o F$ is defined to be the
supremum of $\chi(\o F,b)$, for $b$ ranging over $B$. (This
supremum is either a non-negative real number or $\infty$. Its
value obviously depends on the choice of metric on $A$.) Clearly
$\chi(\o F)=0$ if and only if $\o F$ is continuous.

It will be of interest to have an analog of $\chi(\o F, b)$ for
uniform continuity. If $\o F\FROM B\TO A$ is uniformly continuous,
then for each real $\varepsilon>0$ there exists real $\delta>0$
such that if $U$ is any $\delta$-ball in $B$, then $\o F[U]$ has
diameter $\leq\varepsilon$. We define
\begin{align}             \label{eq-chi-u}
      \chi_u(\o F) \EQ \inf_{\delta>0}\; \sup_{b\in B}\;
               \text{diameter}\; \o F [B_{\delta}(b)],
\end{align}
where the subscript $u$ stands for ``uniform,'' and where
$B_{\delta}(b)$ stands for the $\delta$-ball in $B$ centered at
$b$. This quantity may be called the {\em uniform jump of\/ $\o F$
on $B$.} By analogy with Lemma \ref{lem-continuity}, we have that
$\o F$ is uniformly continuous on $B$ iff $\chi_u(\o F)=0$ (formal
statement and proof omitted). The following lemma is a slight
extension of the well-known equivalence, for compact spaces, of
continuity with uniform continuity. (It may also be well known.)

\begin{lemma}           \label{lem-uniform-compact}
       $\chi(\o F) \,\leq\, \chi_u(\o F)$, with
       equality holding for compact spaces.
\end{lemma}\begin{Proof}
Removing the $b$-supremum from (\ref{eq-chi-u}) we immediately
have, for each $b\in B$,
\begin{align*}
           \chi_u(\o F) &\;\geq \;\inf_{\delta>0}\,\text{diameter}
                     \; \o F [B_{\delta}(b)]\\
           &\EQ \inf \,\{\, \text{diameter}\,\o F[U]\,:\,
                              U \text{ open},\; b\in U \,\}\EQ
                              \chi(\o F,b).
\end{align*}
It follows immediately that
\begin{align*}
            \chi_u(\o F) \;&\geq \;\sup_{b\in B}\, \chi(\o F,b)
                  \EQ \chi(\o F).
\end{align*}

We now consider the opposite inequality $\chi_u(\o F)
 \,\leq\,\chi(\o F)\EQ\sup_{b\in B}\allowbreak\chi(\o F,b)$. If the
 supremum on the right is infinite, the result is immediate. So we
 will assume that it is finite. Let us consider arbitrary real
$K\,>\,\sup_{b\in B}\, \chi(\o F,b)$; it will suffice to prove
that $K\geq\chi_u(\o F)$. We are given that $K\,>\,\chi(\o F,b)$
for every $b\in B$. Thus by the definition (\ref{eq-def-chi}), for
every $b\in B$ there exists an open set $U_b$ with $b\in U$ and
such that $\text{diameter}\;\o F[U_b]<K$.

The sets $U_b$ form an open cover of $A$; by compactness there is
a Lebesgue number for this covering. In other words, there exists
$\delta>0$ such that each $\delta$-ball $B_{\delta}(b)$ is a
subset of $U_{b'}$ for some $b'\in B$. Therefore
$\text{diameter}\,\o F[B_{\delta}(b)]\allowbreak<K$ for each $b\in
B$. Therefore, for this one value of $\delta$,
\begin{align*}
               \sup_{b\in B}\;\text{diameter} \;\o F[B_{\delta}(b)]
                           \;\leq\; K.
\end{align*}
So finally we have
\begin{align*}
               \chi_u(\o F) \EQ \inf_{\delta>0}\,
               \sup_{b\in B}\;\text{diameter} \;\o F[B_{\delta}(b)]
                           \;\leq\; K,
\end{align*}
the desired inequality.
\end{Proof}

Let $(A,d)$ be a metric space, and $\mathbf A =(A,\o{F}_t)_{t\in
T}$ an algebra based on $A$. We define
\begin{align*}
        \chi(\mathbf A,d) \EQ \sup_{t\in T}\; \chi(\o F_t) .
\end{align*}
When the metric $d$ is clear from the context, we may write
$\chi(\mathbf A)$ for  $\chi(\mathbf A,d)$. It should be clear
that $\chi(\mathbf A,d)=0$ if and only if $\mathbf A$ is a
topological algebra.

Finally, for $(A,d)$ a metric space, and $\Sigma$ a set of
equations of similarity type $\la n_t\,:\, t\in T\ra$, we define
\begin{align}                    \label{eq-def-mu}
         \mu(A,d,\Sigma) \EQ \inf\,\{\, \chi(\mathbf A,d) \,:\,
                 \mathbf A = (A,\o{F}_t)_{t\in T}
                    \models\Sigma\;\};
\end{align}
in other words, it is the infimum taken over all algebras built on
$A$ that satisfy $\Sigma$.  When the metric $d$ is clear from the
context, we may write $\mu(A,\Sigma)$ for $\mu(A,d,\Sigma)$.

All these notions have uniform versions, denoted with a subscript
$u$, based on $\chi_u$ in place of $\chi$. Most of this paper
deals with compact metric spaces, on which the two concepts
coincide. In some cases (for example see
\S\S\ref{sub-mult-zero-one}--\ref{sub-majority}) it turns out to
be easier to prove an estimate for $\chi_u$ than for $\chi$. One
should bear in mind that, even when $A$ has been given a metric,
the value of $\chi_u(\mathbf A)$ still depends on the metric that
is chosen to represent the topology on the finite powers
$A^{n_t}$. By Lemma \ref{lem-uniform-compact}, however, there is
no dependence in the compact case. We will mostly work in the
compact case.

In this paper we shall apply (\ref{eq-def-mu}) only in situations
where (i) $A$ is infinite, (ii) $\Sigma$ is finite or countable,
and (iii) $\Sigma$ defines a consistent equational theory (i.e. it
has a model of more than one element). In these circumstanes,
there is at least one model of $\Sigma$ based on $A$; in other
words, the infimum appearing in (\ref{eq-def-mu}) is over a
non-empty set. If one or more of (i--iii) should fail, then it is
possible for the infimum of (\ref{eq-def-mu}) to be over the empty
set. In that case, we would naturally define $\mu(A,d,\Sigma)$ to
be $\infty$.

This quantity $\mu(A,d,\Sigma)$ will be the main object of study
in this paper. Like the previously studied
$\lambda_A(\Sigma)$---see \S\ref{sub-metric-approx-compat}---
$\mu(A,d,\Sigma)$ measures deviation from $A\models\Sigma$, as
follows: in every model of $\Sigma$ based on $A$, some operation
has a discontinuity at least as large as $\mu(A,d,\Sigma)$. (ONLY
APPROXIMATELY TRUE)

It should be apparent that if $A\models\Sigma$, then
$\mu(A,d,\Sigma) = 0$. The converse is false, as we shall see.
Under the right circumstances, there are some connections between
the values of $\lambda_A(\Sigma)$ and $\mu(A,\Sigma)$. See e.g.\
Corollaries \ref{cor-comp-lambda-mu-simple-circle} and
\ref{cor-comp-lambda-mu-k-circle} of \S\ref{sub-close-conseq}.

\section{Elementary remarks about $\mu(A,d,\Sigma)$.}
                                      \label{sec-elem-remarks}
I do not yet have a clear idea of the full scope of
\S\ref{sec-elem-remarks}. Maybe it won't really be necessary, but
for the moment I will file remarks here as I think of them.
\subsection{$\Sigma\subseteq\Sigma'$.}
It is almost obvious that if $\Sigma\subseteq\Sigma'$ then
$\mu(A,d,\Sigma)\leq\mu(A,d,\Sigma')$. (This is simply because, in
evaluating the infimum in (\ref{eq-def-mu}), the set of algebras
for $\Sigma'$ is a subset of the set of algebras for $\Sigma$.)

It is also immediate from (\ref{eq-def-mu}) that if $\Sigma'$ is
any collection of logical consequences of $\Sigma$, then
$\mu(A,d,\Sigma)\geq\mu(A,d,\Sigma')$. Thus, if $\Sigma'$ includes
all of $\Sigma$ together with any subset of the consequences of
$\Sigma$, then $\mu(A,d,\Sigma)=\mu(A,d,\Sigma')$.

It thus follows that, unlike $\lambda_A(\Sigma)$ (see
\S\ref{sub-metric-approx-compat}),  $\mu(A,d,\Sigma)$ is a logical
invariant of $\Sigma$.

\subsection{Topological products.} To come.

\subsection{Products of theories.} To come.

\section{Some sample values of $\mu(A,d,\Sigma)$.}

\subsection{An injective binary operation: $\mu(A,d,\Sigma) = 0$.}
                              \label{subsub-inf-dim}
Consider $\Sigma$ consisting of the two equations
\begin{align}          \label{eq-squared}
                F_0(G(x_0,x_1))\WAVY x_0,\quad\quad
                F_1(G(x_0,x_1)) \WAVY x_1. 
\end{align}
They imply, among other things, that in any topological model $\o
G$ must be a one-one continuous binary operation. Euclidean spaces
of non-zero finite dimension do not have such operations, hence
are not compatible with $\Sigma$. In \S\ref{subsub-inf-dim} we
will be concerned with $A=[0,1]$. Although this $A$ is not
compatible with $\Sigma$, we shall show that
$\mu(A,\delta,\Sigma)=0$ (where $d$ is the ordinary Euclidean
metric).

To begin we let $\o P$ be a continuous function mapping $A=[0,1]$
onto $A^2$. (Such an area-filling curves was devised by G. Peano
in 1890---see \cite[pp.\ 116--7]{wtaylor-gcg}, or many other
sources.) By the Axiom of Choice, $\o P$ has a (discontinuous)
inverse $\o H$; in other words we have
\begin{align*}         
             A^2 \stackrel{\o H}{\TO} A\stackrel{\o P}{\TO} A^2
\end{align*}
with $\o P\Compos\o H$ the identity function on $A^2$.
Let $\lambda=\chi(\o H)$. We remarked above that $\o H$ cannot be
continuous; hence $0<\lambda\leq 1$.

Now to establish our claim that $\mu(A,\Sigma)=0$, we must prove
that the infimum appearing in (\ref{eq-def-mu}) is zero. It will
suffice, given $\varepsilon>0$, to exhibit an algebra $\mathbf A$
based on $A$ with $\chi(\mathbf A)\leq \varepsilon$, and such that
$A\models\Sigma$.

Our algebra $\mathbf A$ is as follows. Its binary operation $\o G$
is defined via
\begin{align*}
               \o G(a_0,a_1) \EQ \varepsilon \,\o H(a_0,a_1).
\end{align*}
The unary operations $\o F_0$ and $\o F_1$ are defined to be the
two components of the function $\o F\FROM A\TO A^2$ that is
defined via
\begin{align*}
               \o F(a) \EQ \o P(1\meet (a/\varepsilon)),
\end{align*}
where $\meet$ denotes the smaller of two real numbers. Clearly $\o
F$ is well-defined and continuous; hence the same is true of $\o
F_0$ and $\o F_1$.

To verify $\Sigma$ for the operations, we first calculate
\begin{align*}
             \o F(\o G(a_0,a_1)) &\EQ
                    \o P(1\meet[\varepsilon \o H(a_0,a_1)/\varepsilon]) \\
              &\EQ  \o P(1\meet[\o H(a_0,a_1)])
                    \EQ  \o P(\o H(a_0,a_1)) \EQ (a_0,a_1),
\end{align*}
with the final equation true by our choice of $\o H$. We now have
$\o F\Compos\o G = 1$, which is tantamount to the equations
$\Sigma$.

As for $\chi(\mathbf A)$, we first note that $\chi(\o G)\leq
\varepsilon\lambda\leq\varepsilon$. By continuity, $\chi(\o F_0) =
\chi(\o F_1) = 0$. Thus $\chi(\mathbf A)\leq\varepsilon$. This
completes the description of this example.

\subsection{$\Sigma$ non-Abelian and simple; $A=S^1$}
                    \label{sub-nonAbel-circle}
Following \cite[\S3.2.3]{wtaylor-asi}, we define a set $\Sigma$ of
equations to be {\em Abelian} iff it is interpretable (in the
sense of \cite{ogwt-mem}) in the equational theory of Abelian
groups. Equivalently, $\Sigma$ is Abelian if and only if it has a
model based on $\Integers$ with operations of the form
\begin{align}          \label{eq-def-op-abelian}
   \o F(x_1,\cdots,x_n)\EQ m_1x_1+\cdots+m_nx_n,
\end{align}
where each $m_i\in\Integers$.

It was proved in Theorem 41 on page 234 of \cite{wtaylor-sae} that
$\Sigma$ is Abelian iff $\Sigma$ is compatible with $S^1$, and
then in \S3.2.3 of \cite{wtaylor-asi}  that $\Sigma$ is Abelian
iff $\lambda_{S^1}(\Sigma)=0$. In fact, if $S^1$ is given its
natural metric as a circle embedded in the Euclidean space
$\Reals^2$, scaled to diameter $1$, then we have
\begin{align}              \label{eq-simple-lambda-1}
      \lambda_{S^1}(\Sigma) \EQ
            \begin{cases}
                  0 & \text{if $\Sigma$ is Abelian}\\
                  1 & \text{otherwise}.
            \end{cases}
\end{align}
The first assertion of (\ref{eq-simple-lambda-1}) obviously holds
for $\mu(S^1,\Sigma)$---namely that $\mu(S^1,\Sigma)=0$ if
$\Sigma$ is Abelian---but the corresponding second assertion is
false.

By a {\em simple term} in the language defined by the operation
symbols $F_t$ of \S\ref{sub-compat-backg}, we mean\footnote%
{This terminology was used, perhaps for the first time, in
Garc\'{\i}a and Taylor \cite{ogwt-mem}, and then again by Taylor
in \cite{wtaylor-eri}.} %
a term that contains at most one $F_t$, and moreover contains at
most one instance of that $F_t$. In other words, according to the
usual recursive definition of terms, a simple term is either a
variable or created at the first stage beyond the inclusion of
variables. An equation $\sigma\wavy\tau$ is {\em simple} iff both
$\sigma$ and $\tau$ are simple terms.

We cannot prove anything like (\ref{eq-simple-lambda-1}) for $\mu$
in place of $\lambda$. However, for some simple non-Abelian
theories $\Sigma$, we can prove the surprisingly exact result that
$\mu(S^1,\Sigma)=2/3$. See
\S\S\ref{sub-mult-zero-one}--\ref{sub-majority}.

\subsection{$\Sigma =$ multiplication with zero and one; $A=S^1$.}
            \label{sub-mult-zero-one}
Throughout \S\ref{sub-mult-zero-one} we let $\Sigma$ be the theory
of a single binary operation with a zero and a one. (Specifically
a left zero and a left one.) Specifically, $\Sigma$ is given by
these equations:
\begin{align}                \label{eq-zero-one}
          F(0,x)\Wavy 0;\quad\quad F(1,x)\Wavy x.
\end{align}
Let us represent the one-sphere $S^1$ as a circle of circumference
2 (hence radius $1/\pi$). We then give it the metric of arc
length: $d(P,Q)$ is the length of the shorter of the two circular
arcs joining $A$ and $B$. In this metric, the space has diameter
$1$. By an {\em arc} in this space we mean the smaller of two
circular arcs joining two points, considered as a closed subset of
$S^1$. We shall prove that, in this metric, $\mu(S^1,\Sigma)=2/3$.

\begin{lemma}                \label{lem-23}
If $F$ is a finite subset of $S^1$ with $\text{diameter}(F)<2/3$,
then there is an arc $A$ of $S^1$, of length
$=\text{diameter}(F)$, such that $F\subseteq A$.
\end{lemma}
\begin{Proof}
Let $\lambda$ denote $\text{diameter}(F)$. Choose $P,Q\in F$ such
that $d(P,Q)=\lambda$, and define $A$ to be the arc
$\overline{PQ}$. Clearly $A$ has length $\lambda$. To prove that
$F\subseteq A$, we consider three intervals of length $\lambda$.
The first is $A$ itself; the second is $B$ which meets $A$ only in
the point $P$; the third is $C$ which meets $A$ only at the point
$Q$. Since $\lambda<2/3$, the intervals $B$ and $C$ are disjoint.
Now every member of $F$is within $\lambda$ of $P$, hence belongs
either to $A$ or to $B$. Similarly every member of $F$ belongs
either to $A$ or to $C$. Now suppose that $R\in F$ but $\not\in
A$. Then $R$ must belong to both $B$ and $C$, which is a
contradiction; this contradiction completes the proof that
$F\subseteq A$.
\end{Proof}

{\bf Remarks.} The proof can easily be extended to infinite $F$,
although we will not need this refinement. The number $2/3$ is
sharp for this lemma, as follows. Consider $F$ whose members are
three points equally spaced at distance $2/3$ about the circle
$S^1$. Clearly no arc contains $F$, but $\text{diameter}\,(F) =
2/3$.

\begin{theorem}           \label{thm-23}
 $\mu(S^1,\Sigma)\EQ 2/3$. (With $\Sigma$ as defined in
 (\ref{eq-zero-one}).)
\end{theorem}\begin{Proof}
{\bf Part 1.} $\mu(S^1,\Sigma)\leq 2/3$. We must exhibit an
algebra $\mathbf A=\la S^1,\o F,\o0,\o1\ra$ (with $\o F$ binary)
satisfying equations (\ref{eq-zero-one}), and with $\chi(\o
F)\leq2/3$. For convenience, let us take $R$ to stand for $1/\pi$,
the radius of our circle $S^1$. We now define the operations of
$\mathbf A$ as follows: $\o1= R$, $\;\o0=-R$, and $\o F$ is given
by these formulas:
\begin{align}
           \o F(R,z) \EQ z; \quad&\quad \o F(-R,z)\EQ -R;
                     \label{eq-def-0-1-S1-a}\\
           \o F(w\neq\pm R,\, Re^{i\theta})&\EQ
 \begin{cases}
          R\, e^{i\pi/3} & \text{if}\;\;0\leq\theta<2\pi/3,\\
          -R           & \text{if}\;\;2\pi/3\leq\theta<4\pi/3,\\
          R\,e^{5i\pi/3}  & \text{if}\;\,4\pi/3\leq\theta<6\pi/3.\\
         \end{cases}     \label{eq-def-0-1-S1-b}
\end{align}
The satisfaction of (\ref{eq-zero-one}) is immediate from
(\ref{eq-def-0-1-S1-a}) of the definition.

To evaluate $\chi(\o F)$, we first consider $\chi(\o F, (w,z))$,
where $w\neq R$. It is obvious from
(\ref{eq-def-0-1-S1-a}--\ref{eq-def-0-1-S1-b}) that $(w,z)$ has a
neighborhood $U$ such that $\o F[U]\subseteq R\{-1, e^{i\pi/3},
e^{5i\pi/3}\}$. This latter set has diameter $2/3$, and so we may
turn our attention to the case of $w=R$.

We begin the case of $w=R$ by remarking that $\o F$ has a 3-fold
symmetry, as follows: if $w\neq -R$, then $\o F(w,e^{2\pi
i/3}z)\,=\,e^{2\pi i/3}\o F(w,z)$. Since multiplication by any
unimodular complex scalar is a rotation, and does not change
diameters, it will be enough to focus our attention on $\o
F(R,Re^{i\theta})$ where $0\leq \theta<2\pi i/3$. We first assume
that $\theta\neq 0$. We may then consider a neighborhood
$U_0\times U_1$ of $(R,Re^{i\theta})$, where $-R\not\in U_0$ , and
$U_1$ is a small arc about $Re^{i\theta}$ that lies interior to
the arc $\o{R\,Re^{2\pi i/3}}$. Then $\o F [U_0\times U_1]$ is
$U_1 \cup \{Re^{\pi i/3}\}$. The reader may easily check that this
set has diameter $<2/3$.

It finally remains to consider $w=R$ and $\theta=0$, which is to
say, to evaluate $\chi(\o F,(R,R))$. Things go exactly as before,
except that a neighborhood $U_0\times U_1$ of $(R,R)$ will contain
some points of the form $(v,Re^{i\theta})$ where $v\neq\pm R$ and
$4\pi/3<\theta<2\pi$. Thus $\o F [U_0\times U_1]$ will be $U_1
\cup \{Re^{\pi i/3},\, Re^{5\pi i/3}\}$.\vspace{0.1in}

{\bf Part 2.} $\mu(S^1,\Sigma)\geq 2/3$.  For a proof by
contradiction, let us assume that $\mu(S^1,\Sigma)<2/3$. By
(\ref{eq-def-mu}), there is an algebra $\mathbf A=\la S^1,\o
F,\o0,\o1\ra$ such that $\chi(\o F)<2/3$ and such that $\mathbf
A\models\Sigma$. Let us give $(S^1)^2$ the sum metric
\begin{align}            \label{eq-sum-metric}
d((a,b),\,(c,d)) \EQ d(a,c)\,+\,d(b,d).
\end{align}
 By Lemma \ref{lem-uniform-compact}, we also have $\chi_u(\o F)<2/3$.
Referring to (\ref{eq-chi-u}) (the definition of $\chi_u$), we see
that there exists $\delta>0$ such that $d((a,b),(c,d))<\delta$
implies $d(\o F(a,b),\o F(c,d))<2/3$. Let
\begin{align*}
             t_0, \; t_1, \; \cdots \;t_{N-1}, \; t_N\EQ t_0
\end{align*}
be points of $S^1$ such that
\begin{itemize}
\item[(a)] The points $t_i$ are evenly spaced around the circle,
   with $d(t_i,t_{i+1})\,<\,\delta/2$ for all appropriate $i$. We
   will refer to the portion of $S^1$ between $t_i$ and $t_{i+1}$
   as a {\em segment} of the circle.
\item[(b)] This sequence of points continues around the circle in
the same direction, and goes around the circle exactly once.
\item[(c)] For convenience, we make sure that $\o0= t_0$ and $\o1=
t_K$ for some $K$.
\end{itemize}

For the remainder of the proof we consider the restriction of $\o
F$ to the finite set $\{(t_i,t_j)\,:\, 0\leq i\leq K,\; 0\leq j<
N\,\}$. From (a) and (\ref{eq-sum-metric}) and our choice of
$\delta$, we immediately have
\begin{align}                \label{eq-nearby}
      \text{diameter}\; \{\o F(t_i,t_j), \o F(t_i,t_{j+1}), \o
      F(t_{i+1},t_j), \o F(t_{i+1},t_{j+1})\,\} \;<\; 2/3.
\end{align}
for $0\leq i<K,\; 0\leq j< N$. We will finish the proof by showing
that the metric arrangement (\ref{eq-nearby}) and Equations
(\ref{eq-zero-one}) together lead to a contradiction.

Applying Lemma \ref{lem-23} to (\ref{eq-nearby}) we see that for
$0\leq i<K,\; 0\leq j< N$ there is an arc $A_{ij}$ of length
$<2/3$ such that
\begin{align}                  \label{eq-near-in-arc}
    \o F(t_i,t_j),\, \o F(t_i,t_{j+1}), \,\o
      F(t_{i+1},t_j), \,\o F(t_{i+1},t_{j+1}) \;\in\; A_{ij}.
\end{align}

For $i=0,\ldots,\, K$, let us take a continuous function
$\o\gamma_i\FROM S^1\TO S^1$, in such a way that the following
five conditions are met:
\begin{itemize}
\item[(i)] $\o\gamma_i(t_j) = \o F(t_i,t_j)$
           and $\o\gamma_i(t_{j+1}) = \o F(t_i,t_{j+1})$;
\item[(ii)] For $0\leq i <K$, $\;\o\gamma_i$ maps the arc $\o{t_jt_{j+1}}$
       into the arc $A_{ij}$.
\item[(iii)] For $0<i\leq K$, $\;\o\gamma_i$ maps the arc $\o{t_jt_{j+1}}$
       into the arc $A_{(i-1)j}$.
\item[(iv)] $\o\gamma_0$ is the constantly $\o0$ function.
\item[(v)] $\o\gamma_K$ is the identity function.
\end{itemize}
(Condition (i) can be met directly. For conditions (ii) and (iii),
we use (\ref{eq-near-in-arc}) to see that the endpoints $t_j$ and
$t_{j+1}$ both map into the arc $A_{ij}\cap A_{(i-1)j}$; hence the
arc between them can be mapped into $A_{ij}\cap A_{(i-1)j}$ (here
we use the fact that a non-empty intersection of two arcs of
diameter $<1$ is itself an arc). For condition (iv), we recall
that $\o0$ is $t_0$; therefore the equations (\ref{eq-zero-one}),
together with condition (i), tell us that all the values
$\o\gamma_0(t_j)$ are $\o0$. Therefore we easily satisfy
conditions (ii) and (iii) by making $\gamma_0$ constantly equal to
$\o0$---which yields also (iv). Condition (v) is satisfied
similarly.)

Now, for a contradiction, we will prove that $\o\gamma_0$ is
homotopic to $\o\gamma_K$, in contradiction to conditions (iv) and
(v). Using the transitivity of homotopy, it will be enough to
prove that $\o\gamma_i$ is homotopic to $\o\gamma_{i+1}$ for
$0\leq i<K$. So we fix a value of $i$ in this range, and proceed
to define the required homotopy.

For $0\leq j<N$, we define a continuous $S^1$-valued function $\o
G_j$, whose domain is the arc $\o{t_it_{i+1}}$ and which satisfies
\begin{itemize}
 \item[(vi)] $\o G_j(t_i)=\o
          F(t_i,t_j)$ and $\o G_j(t_{i+1})=\o F(t_{i+1},t_j)$.
 \item[(vii)] $\o G_j$ maps the arc $\o{t_it_{i+1}}$ into the arc
           $A_{ij}$.
 \item[(viii)] $\o G_{j+1}$ maps the arc $\o{t_it_{i+1}}$ into the arc
           $A_{ij}$.
\end{itemize}
(Again, these conditions are all possible by
(\ref{eq-near-in-arc}).)

We now consider the set $B_{ij}\,=\,\o{t_it_{i+1}}
\times\o{t_jt_{j+1}}\,\subseteq\, S^1\times S^1$. We define a
function $\o\phi_{ij}$ from the boundary of $B_{ij}$ to the arc
$A_{ij}$, as follows:
\begin{align}
\o\phi_{ij}(s,t_j)&\EQ \o G_j(s) \label{phiij-one} \\
 \o\phi_{ij}(s,t_{j+1})&\EQ \o G_{j+1}(s)\label{phiij-two} \\
\o\phi_{ij}(t_i,t)&\EQ\o\gamma_i(t)\label{phiij-three} \\
 \o\phi_{ij}(t_{i+1},t)&\EQ\o\gamma_{i+1}(t)\label{phiij-four}
\end{align}
(The reader may check, from what has come before, that
$\text{Range}\,(\o\phi_{ij})\,\subseteq \,A_{ij}$, and that
$\o\phi_{ij}$ is well-defined, and hence continuous, at the
corners of $B_{ij}$.)

Any continuous function from the boundary of a plane disk to the
real line extends to a continuous function defined on the full
disk. (This is Tietze's Extension Theorem.) Thus there exists a
continuous function $\o\Phi_{ij}\FROM B_{ij}\TO A_{ij}$ that
restricts to $\o\phi_{ij}$ on the boundary.

We will now show that
\begin{align*}
               \o\Phi_i \EQ \bigcup_{j=0}^{N-1} \o\Phi_{ij}
\end{align*}
is the desired homotopy between $\o\gamma_i$ and $\o\gamma_{i+1}$.
Clearly its domain is $\bigcup_{j=0}^{N-1} B_{ij} \;\allowbreak=\;
\o{t_it_{i+1}}\times S^1$, and by
(\ref{phiij-one}--\ref{phiij-two}), for each $j$ the component
functions $\o\phi_{ij}$ and  $\o\phi_{i(j+1)}$ agree where they
overlap. Thus $\o\Phi_i$ is a continuous function defined on
$\o{t_it_{i+1}}\times S^1$. Finally, from
(\ref{phiij-three}--\ref{phiij-four}) it follows that, for all
$t\in S^1$, we have $\o\Phi_i(t_i,t)=\gamma_i(t)$ and
$\o\Phi_i(t_{i+1},t)=\gamma_{i+1}(t)$. Thus $\o\Phi_i$ is the
desired homotopy. As mentioned above, transitivity yields a
homotopy between the identity and a constant function. This
contradiction to known results completes the proof of the theorem.
\end{Proof}

\hspace{-\parindent}%
{\bf Remark on the proof.} In fact, what we have done here
is---for a contradic\-tion---to begin with a solution $\o F$ to
the equations $\Sigma$, such that $\o F$ is discontinuous, but by
no more than $2/3$. We have then focused on a finite subset of $\o
F$ (comprising the function-values $\o F(t_i,t_j)$). Finally we
have interpolated a continuous function $\o G$ through these
values that also satisfies $\Sigma$. Since no such $\o G$ exists,
we have our contradiction. In \S\ref{sub-idem-comm} and
\S\ref{sub-majority} we will see this method to be widely
applicable. See also \S\ref{sub-proofs-comment}.

\subsection{$\Sigma=$ commutative idempotent binary; $A=S^1$}
                   \label{sub-idem-comm}
In \S\ref{sub-idem-comm} we consider the following $\Sigma$, which
defines commutative idempotent binary operations:
\begin{align}                   \label{eq-idem-comm}
        F(x,y) \Wavy F(y,x); \quad\quad F(x,x)\Wavy x.
\end{align}
In \S3.4.1 of \cite{wtaylor-asi} we remarked that this $\Sigma$ is
non-Abelian, hence not compatible with $S^1$. Since it is also
simple, its $\mu$-value is amenable to estimation. By a method
similar to that of \S\ref{sub-mult-zero-one} (again using Lemma
\ref{lem-23}) we will prove
\begin{theorem}       \label{th-idem-comm-23}
   $\mu(S^1,\Sigma)\EQ2/3$. (With\/ $\Sigma$ as defined in
 (\ref{eq-idem-comm}).)
\end{theorem}
\begin{Proof} \vspace{0.1in}

{\bf Part 1.} $\mu(S^1,\Sigma)\leq 2/3$. We must exhibit an
algebra $\mathbf A=\la S^1,\o F\ra$ (with $\o F$ binary)
satisfying equations (\ref{eq-idem-comm}), and with $\chi(\o
F)\leq2/3$. To avoid fractions, we will represent elements of our
circle as real numbers modulo $3$, and will assume that these
numbers parametrize the distance. In this reframing, the circle
has diameter $3/2$, and so we expect to prove that $\chi(\o F)\leq
1$.
 We define the operation $\o F$ of
$\mathbf A$ as follows:
\begin{align}                \label{def-F-idem-comm}
              \o F(s,t)\EQ \o F(t,s) \EQ\begin{cases}
                     s\join t &\text{if $\;0\leq s,t\leq 1$} \\
                       s\join t  &\text{if $\;1\leq s,t\leq3$} \\
                    1\,+\,s\,+\,t
                                  &\text{if $\;0\leq s\leq 1$ and
                                           $2\leq t\leq3$}\\
                     2+t  &\text{if $\;0\leq s\leq 1$ and
                                           $1\leq t\leq2$,}
              \end{cases}
\end{align}
where of course the addition is taken modulo $3$. It is obvious
that this $\o F$ satisfies the $\Sigma$ in (\ref{eq-idem-comm}).
In order to estimate $\chi(\o F)$ we consider the following
diagram:\vspace{0.2in}

\begin{align}        \label{eq-sudoku}
\begin{tabular}{|l|l|l|}
\hline \begin{minipage}[t]{2cm}1\hfill 2\\[1.2cm]
          0\hfill 1\end{minipage} &
           \begin{minipage}[t]{2cm}3\hfill 3\\[1.2cm]
          2\hfill 2\end{minipage}&
           \begin{minipage}[t]{2cm}3\hfill 3\\[1.2cm]
          2\hfill 3\end{minipage}\\
\hline \begin{minipage}[t]{2cm}1\hfill 1\\[1.2cm]
          0\hfill 0\end{minipage} &
           \begin{minipage}[t]{2cm}2\hfill 2\\[1.2cm]
          1\hfill 2\end{minipage}&
           \begin{minipage}[t]{2cm}2\hfill 3\\[1.2cm]
          2\hfill 3\end{minipage}\\\hline \begin{minipage}[t]{2cm}1\hfill 1\\[1.2cm]
          0\hfill 1\end{minipage} &
           \begin{minipage}[t]{2cm}0\hfill 1\\[1.2cm]
          0\hfill 1\end{minipage}&
           \begin{minipage}[t]{2cm}1\hfill 2\\[1.2cm]
          0\hfill 1\end{minipage}\\
          \hline
\end{tabular} \vspace{0.2in}
\end{align}

This illustration depicts $[0,3]\times[0,3]$, divided into nine
squares of dimensions $1\times 1$. If we consider the diagram
modulo $3$ in each direction, then we have our version of the
torus $S^1\times S^1$. As the reader may check---case by
case---the four values shown in each small square are the corner
$\o F$-values for that square, as supplied by our definition
(\ref{def-F-idem-comm}). Moreover, on each small edge, the $\o
F$-values (considered not modulo $3$, but as reals in $[0,3]$)
vary linearly between the indicated corner values.

It is now not hard to observe---again, case by case---that no jump
is greater than $1$ in the limit. The most serious case occurs at
the upper-right corner, call it $P$, of the upper-left small
square: the values at $P$ are $0$, $1$, $2$, $3=0$. Given
$\varepsilon>0$, there is a neighborhood $U$ of $P$ such that $\o
F[U]\subseteq [-\varepsilon,\varepsilon]\cup [1-\varepsilon,
1+\varepsilon] \cup [2-\varepsilon, 2+\varepsilon]$. This last set
has diameter $1+2\varepsilon$. Then $\chi(\o F,P)$ is the relevant
infimum, which clearly is $1$. We have now established the
required properties of $\mathbf A=\la S^1,\o F\ra$, and hence the
proof of Part 1 is complete.\footnote{%
Actually, formulas (\ref{def-F-idem-comm}) are not especially
relevant or important to the proof. The important thing is the
sudoku-like puzzle of finding Diagram (\ref{eq-sudoku}): the
values shown must illustrate idempotence, commutativity and small
jumps. From there one can easily contrive a function like our $\o
F$.} \vspace{0.1in}

{\bf Part 2.}  $\mu(S^1,\Sigma)\geq 2/3$ For a proof by
contradiction, we assume that $\mu(S^1,\Sigma)\;<\;2/3$.  By
(\ref{eq-def-mu}), there is an algebra $\mathbf A=\la S^1,\o F\ra$
such that $\chi(\o F)<2/3$ and such that $\mathbf A\models\Sigma$.
As in the proof of Theorem \ref{thm-23}, we give $(S^1)^2$ the sum
metric. As before, there exists $\delta>0$ such that
$d((a,b),(c,d))<\delta$ implies $d(\o F(a,b),\o F(c,d))<2/3$. Let
\begin{align*}
             t_0, \; t_1, \; \cdots \;t_{N-1}, \; t_N\EQ t_0
\end{align*}
be points of $S^1$ satisfying (a--c)  in the proof of Theorem
\ref{thm-23}. As before, we have\footnote{%
The ``+1'' appearing in subscripts in
(\ref{eq-nearby-repeat}--\ref{eq-near-in-arc-repeat}), and
elsewhere, is of course to be understood modulo $N$.}
\begin{align}                \label{eq-nearby-repeat}
      \text{diameter}\; \{\o F(t_i,t_j), \o F(t_i,t_{j+1}), \o
      F(t_{i+1},t_j), \o F(t_{i+1},t_{j+1})\,\} \;<\; 2/3.
\end{align}
for $0\leq i, j< N$. The nearness relations
(\ref{eq-nearby-repeat}) will make it possible to define a
continuous binary operation $\o G$ that interpolates the $N^2$
discrete function values $\o F(t_i,t_j)$ ($0\leq i,j<N$). Using
the fact that these values obey (\ref{eq-idem-comm}) we will be
able to make sure that the interpolated operation $\o G$ also
obeys (\ref{eq-idem-comm}). Thus we will have $\la S^1,\o
G\ra\models \Sigma$ with $\o G$ continuous, in contradiction to
the known fact \cite[\S3.4.1]{wtaylor-asi} that
$S^1\not\models\Sigma$; this contradiction will complete the proof
of the theorem.

Applying Lemma \ref{lem-23} to (\ref{eq-nearby-repeat}) we see
that for $0\leq i, j< N$ there is an arc $A_{ij}$ of length $<2/3$
such that
\begin{align}                  \label{eq-near-in-arc-repeat}
    \o F(t_i,t_j),\, \o F(t_i,t_{j+1}), \,\o
      F(t_{i+1},t_j), \,\o F(t_{i+1},t_{j+1}) \;\in\; A_{ij}.
\end{align}
From (\ref{eq-near-in-arc-repeat}) we easily derive, for $0\leq i,
j< N$, that
\begin{align}
           \o F(t_i,t_j),\, \o F(t_i,t_{j+1}) \;&\in\;
                            A_{ij}\,\cap\,A_{(i-1)j}
                               \label{eq-limit-horiz} \\
           \o F(t_i,t_j),\, \o F(t_{i+1},t_j) \;&\in\;
                            A_{ij}\,\cap\,A_{i(j-1)}.
                                 \label{eq-limit-vert}
\end{align}
From (\ref{eq-near-in-arc-repeat}) and Equations
(\ref{eq-idem-comm}) we also have
\begin{align}             \label{eq-idempot-on-aij}
             t_i \;\in\; A_{ij}\cap A_{i(j-1)} \cap A_{(i-1)j}
                 \cap A_{(i-1)(j-1)}
\end{align}
for all $i$. By symmetry (\ref{eq-idem-comm}), we have $\o
F(t_i,t_j)\,=\,\o F(t_j,t_i)$ for all $i$ and $j$; hence we may
further require that
\begin{align}       \label{eq-symm-Aij}
        A_{ij} \EQ A_{ji}
\end{align}
for all appropriate $i$ and $j$. Moreover, since each $A_{ij}$ is
an arc of length $\leq2/3$, the right-hand sides of
(\ref{eq-limit-horiz}), (\ref{eq-limit-vert}) and
(\ref{eq-idempot-on-aij}) are themselves arcs of $S^1$.

Turning to the definition of $\o G$, we begin with what may be
called the coordinate circles, $(s,t_j)$ and $(t_i,t)$ for $s,t\in
S^1$ and $0\leq i,j < N$. From (\ref{eq-limit-horiz}) and
(\ref{eq-limit-vert}) it is clear that for these domain values we
may now define a continuous binary operation $\o G$ satisfying the
following conditions for all $s,t\in S^1$ and $0\leq i,j < N$:
\begin{itemize}
 \item[(i)] $\o G(t_i,t_j) \EQ \o F(t_i,t_j)$;
 \item[(ii)] $\o G(s,t_j) \in\;A_{ij}\,\cap\,A_{i(j-1)}\,,\;$
             for $s$ in the arc $\o{t_i\,t_{i+1}}\;$;
 \item[(iii)] $\o G(t_i,t) \in\;A_{ij}\,\cap\,A_{(i-1)j},\;$
        for $t$ in the arc $\o{t_j\,t_{j+1}}\;$.
\end{itemize}
By the symmetry that we already have, e.g. (\ref{eq-symm-Aij}), we
may further require
\begin{itemize}
 \item[(iv)] $\o G(s,t_j) \EQ \o G(t_j,s)$
\end{itemize}
for all appropriate $j$ and $s$.

For each $i,j$ we have defined $\o G$ on the boundary of the
rectangle $\o{t_i\,t_{i+1}}\times \o{t_j\,t_{j+1}}$, which
consists of the following four arcs:
\begin{align}
   \{t_i\}\!\times\!\o{t_j\,t_{j+1}}\,,\;\;\;\;
   \o{t_i\,t_{i+1}}\!\times\!\{t_{j+1}\}\,,\;\;\;\;
    \{t_{i+1}\}\!\times\!\o{t_j\,t_{j+1}}\,,\;\;\;\;
     \o{t_i\,t_{i+1}}\!\times\!\{t_{j}\}\,.
\end{align}
By (ii) and (iii), our partial operation $\o G$ maps each of these
four arcs into $A_{ij}$. In other words $\o G$ maps the boundary
of the rectangle $\o{t_i\,t_{i+1}}\times \o{t_j\,t_{j+1}}$ into
the topological interval $A_{ij}$. As is well known, $\o G$ may be
extended to a continuous function on the entire rectangle:
\begin{align*}
             \o G_{ij}\,\FROM\, \o{t_i\,t_{i+1}}\times \o{t_j\,t_{j+1}}
              \TO A_{ij}.
\end{align*}
Let us take such a $\o G_{ij}$ for every $i$ and $j$ with $i\leq
j$. Then for $i>j$ we will define $\o G_{ij}$ by the formula
\begin{align}             \label{eq-patchwork-symm}
            \o G_{ij}(s,t) \EQ \o G_{ji}(t,s).
\end{align}
It is obvious from (iv) that the $\o G_{ij}$ defined by
(\ref{eq-patchwork-symm}) also extends our given $\o G$ as defined
on the boundary of $\o{t_i\,t_{i+1}}\times \o{t_j\,t_{j+1}}$.

It should now be clear that $\bigcup_{0\leq i,j<N}\o G_{ij}$ is a
continuous binary operation on $S^1$ that extends our partial
operation $\o G$. We will denote this full operation also by $\o
G$. For (\ref{eq-idem-comm}), we need to check  its idempotence
and its symmetry. For this, we need to make two further
stipulations in the definition of $\o G_{ii}$ (for $0\leq i<N$).
Since $\o F$ satisfies (\ref{eq-idem-comm}), we have $\o
F(t_i,t_i)\,=\,t_i$ for all $i$. By (i) we have $\o G(t_i,
t_i)\,=\,t_i$ for all $i$. Let us first extend $\o G_{ii}(s,s)$ to
have the value $s$ for each $s\in \o{t_i\,t_{i+1}}$. The diagonal
$\{(s,s):s\in\o{t_i\,t_{i+1}}\;\}$ divides $\o{t_i\,t_{i+1}}\times
\o{t_i\,t_{i+1}}$ into two triangles, and $\o G_{ii}$ has been
defined on the boundary of each of these triangles. Then  $\o
G_{ii}$ may be extended to one triangle (as before), and reflected
to the other triangle by the formula $\o G_{ii}(t,s)\,=\,\o
G_{ii}(t,s)$. This completes a definition of $\o G_{ii}$ on the
full rectangle $\o{t_i\,t_{i+1}}\times \o{t_i\,t_{i+1}}$.

It is now obvious that $\o G$ satisfies (\ref{def-F-idem-comm}) if
the variables are assigned values in any rectangle
$\o{t_i\,t_{i+1}}\times \o{t_i\,t_{i+1}}$. For values outside such
a rectangle, idempotence is moot, and (\ref{eq-patchwork-symm})
suffices to prove symmetry. We have thus constructed a continuous
commutative idempotent operation on $S^1$, in contradiction to
known results. This contradiction completes the proof of the
theorem.
\end{Proof}

\subsection{$\Sigma=$ ternary majority laws;
$A=S^1$} \label{sub-majority}

 In \S\ref{sub-majority} we follow the general path of
\S\ref{sub-mult-zero-one} and \S\ref{sub-idem-comm}, but this time
we consider a non-Abelian simple theory about a {\em ternary}
operation symbol $F$. In this section we let $\Sigma$ denote the
following equations, known sometimes as the {\em majority
equations}:
\begin{align}          \label{eq-majority}
      F(x,x,z) \Wavy F(x,z,x) \Wavy F(z,x,x) \Wavy x.
\end{align}
For the sake of completeness we also consider the {\em symmetric
majority equations:}
\begin{align}          \label{eq-majority-symm}
     \Sigma' \EQ \Sigma \,\cup\, \{F(x,y,z)\wavy F(x,z,y)\wavy
                     F(y,z,x)\}.
\end{align}

 Again using Lemma
\ref{lem-23}, we will prove
\begin{theorem}           \label{th-s1-majority}
$\mu(S^1,\Sigma)\EQ\mu(S^1,\Sigma')\EQ2/3$. (With\/ $\Sigma$,
$\Sigma'$ as defined in (\ref{eq-majority}),
(\ref{eq-majority-symm}), resp.)
\end{theorem}
\begin{Sketch} It will of course be enough to prove that
$2/3\leq\mu(S^1,\Sigma)\leq\mu(S^1,\Sigma')\leq2/3$.
\vspace{0.1in}

{\bf Part 1.} $\mu(S^1,\Sigma')\leq 2/3$. We must exhibit an
algebra $\mathbf A=\la S^1,\o F\ra$ (with $\o F$ ternary)
satisfying equations (\ref{eq-majority-symm}), and with $\chi(\o
F)\leq2/3$. For convenience, as in the proof of Theorem
\ref{thm-23}, we represent $S^1$ as a circle of radius $R \,=\,
1/\pi$, with metric determined by arc length around the circle.
Thus in this representation $S^1$ has diameter $1$.

On $S^1$ we shall construct a ternary operation $\o F$ satisfying
three properties, which guarantee (\ref{eq-majority-symm}) and
which allow us to make the desired estimate of $\chi(\o F)$:
\begin{itemize}
 \item[(i)] $\o F$ satisfies $F(x,y,z)\wavy F(x,z,y)\wavy F(y,z,x)$.
 \item[(ii)] $\o F(a,b,c)\in\{a,b,c\}$ for all $a,b,c$.
 \item[(iii)] If $d(a,a')<2/3$, then $\o F(a,a',b)
         \in \o
           {a\,a'}$.
\end{itemize}



Here is the definition of $\o F$. Given  $a,b,c\in S^1$, we
examine the three distances $d(a,b)$, $d(b,c)$ and $d(c,a)$.
\vspace{0.1in}

{\bf Definition of $\o F$, clause (1).} If all three distances are
$<2/3$, then one of these distances is the sum of the other two.
For example $d(a,c)=d(a,b)+d(b,c)$. In that case $b$ is said to be
between $a$ and $c$, and we define $\o F(a,b,c)$ to be $b$. The
same formula, {\em mutatis mutandis}, yields $a$ (resp.\ $c$)
between the other two, in which case $\o F(a,b,c)$ is $a$ (resp.\
$c$). \vspace{0.1in}

{\bf Definition of $\o F$, clause (2).} If exactly two of the
three distances are $<2/3$, then $a$, $b$ and $c$ must be
distinct, as the reader may verify. For example we might have
$0<d(a,b), d(b,c) < 2/3$ and $d(a,c)>2/3$. In this case, $a$ and
$c$ must lie on opposite sides of $b$, for otherwise $d(a,c)$
would be too small. In this case, we define $\o F(a,b,c)$ to be
$b$. We extend the definition, {\em mutatis mutandis}, to the
other two possible arrangements. \vspace{0.1in}

{\bf Definition of $\o F$, clause (3).} If exactly one of the
three distances is $<2/3$, say $d(a,b)<2/3$, then we define $\o
F(a,b,c)$ to be either $a$ or $b$, chosen at random. We extend the
definition, {\em mutatis mutandis}, to the other two possible
arrangements. \vspace{0.1in}

{\bf Definition of $\o F$, clause (4).} Finally, if none of the
three distances is $<2/3$, then all three must be equal to $2/3$.
In this case we let $\o F(a,b,c)$ be $a$ or $b$ or $c$, chosen at
random. \vspace{0.1in}

We now turn to the verification of (i), (ii) and (iii) for our
operation $\o F$. Condition (i) is immediate, since in all cases
the definition concerns e.g. a set of distances; it does matter in
which order the three variables enter the triple $(a,b,c)$.
Condition (ii) is immediate from the construction of $\o F$.

As for Condition (iii), let us consider the definition of $\o
F(a,a',b)$, where $d(a,a')<2/3$. If $\o F(a,a',b)$ falls into
clause (1) of the definition, then we may discern two cases: (a)
$b$ is between $a$ and $a'$, and (b) it is not. In  case (a), $\o
F(a,a',b)$ is $b$, which lies in the interval $\o{a\,a'}$. In case
(b), $\o F(a,a',b)$ is either $a$ or $a'$, and both of these lie
in the arc $\o{a\,a'}$.

Verifying Condition (iii) for clause (2) of the definition, if
$d(a,a')<2/3$ then we cannot have  $a$ and $a'$ on opposite sides
of $b$ (for then all three interals would be small). Thus either
$a$ and $b$ are on opposite sides of $a'$, or $a'$ and $b$ are on
opposite sides of $a$. Thus we have $\o F(a,a',b)$ equal to $a$ or
$a'$, and hence in the arc $\o{a\,a'}$.

The verification of Condition (iii) for clause (3) of the
definition is immediate. Clause (4) cannot occur in the
calculation of $\o F(a,a',b)$. Hence we have considered all
clauses for the evaluation of $\o F(a,a',b)$; hence Condition
(iii) is verified.

Having established conditions (i--iii), we turn now to our
previous claim that these conditions imply the desired properties
for $\o F$. As for equations (\ref{eq-majority-symm}), Condition
(i) is symmetry itself, and condition (iii) immediately yields the
majority laws  (\ref{eq-majority}). All that remains for Part 1 of
the proof is to estimate $\chi(\o F,(a,b,c))$ for $(a,b,c)\in
(S^1)^3$. Our estimate will be based solely on conditions
(i--iii). We consider two possibilities for the triple $(a,b,c)$.
\vspace{0.1in}

\hspace*{-\parindent}%
{\bf Case 1: $d(a,b)=d(b,c)=d(c,a)=2/3$.} In this case,
$\{a,\,b,\,c\}$ is an equilateral triangle of diameter $2/3$. For
a neighborhood of $(a,b,c)$ in $(S^1)^3$, we may consider a set
$U\times V\times W$, where $U$ (resp.\ $V$, $W$) is a neighborhood
of $a$ (resp.\ $b$, $c$). From Condition (ii), we easily see that
\begin{align*}
              \o F[U\times V\times W] \;\subseteq\;
                   U\,\cup\, V\,\cup\, W.
\end{align*}
By making the neighborhoods $U$, $V$ and $W$ small, we obviously
have $\allowbreak\text{diameter}\allowbreak\,\o
F[U\allowbreak\times V\times W] < 2/3 + \varepsilon$ for any
$\varepsilon>0$. Thus $\chi(\o
F,(a,b,c))\,\leq\,2/3$.\vspace{0.1in}

\hspace*{-\parindent}%
{\bf Case 2: either $d(a,b)\neq 2/3$ or $d(b,c)\neq2/3$ or
$d(c,a)\neq2/3$.} Then obviously one of these three distances must
be $<2/3$. Since $\o F$ is symmetric, we assume without loss of
generality that $d(a,b)<2/3$.  Choose real $\delta$ with
$0<2\delta<(2/3 -d(a,b))$, and let $U$ (resp. $V$) be the
$\delta$-ball about $a$ (resp. $b$) with radius $\delta$. If $u\in
U$ and $v\in V$, then $d(u,v)<2/3$. Hence, for any $w$, we have
$\o F(u,v,w)\in\o{u\,v}$, by (iii). In other words, we have
\begin{align}
            \o F\,[U\times V\times S^1] \;\subseteq\; U\,\cup\,
            V\,\cup\,\o{a\,b}.
\end{align}
This last is a set of diameter $d(a,b)+2\delta$; by our choice of
$\delta$ this diameter $< 2/3$. In other words, we have now shown
that $\chi(\o F,(a,b,c))\,<\,2/3$.

Combining Cases 1 and 2, we see that $\chi(\o F)\leq2/3$, and
hence that $\mu(S^1,\Sigma)\leq2/3$. This finishes Part 1 of the
proof.\vspace{0.1in}

{\bf Part 2.} $\mu(S^1,\Sigma)\geq2/3$. For a proof by
contradiction, we assume that $\mu(S^1,\Sigma)\,<\,2/3$.  By
(\ref{eq-def-mu}), there is an algebra $\mathbf A=\la S^1,\o F\ra$
such that $\chi(\o F)<2/3$ and such that $\mathbf A\models\Sigma$.
As in the proof of Theorems \ref{thm-23} and
\ref{th-idem-comm-23}, we give $(S^1)^3$ the sum metric (in this
case, the sum of distances over three coordinates). As before,
there exists $\delta>0$ such that $d((a,b,c),(d,e,f))<\delta$
implies $d(\o F(a,b,c),\o F(d,e,f))<2/3$. Let
\begin{align*}
             t_0, \; t_1, \; \cdots \;t_{N-1}, \; t_N\EQ t_0
\end{align*}
be points of $S^1$ satisfying (a--c)  in the proof of Theorem
\ref{thm-23}. As before, we have\footnote{%
The ``+1'' appearing in
(\ref{eq-nearby-third}--\ref{eq-near-in-arc-third}), and
elsewhere, is again modulo $N$.}
\begin{align}                \label{eq-nearby-third}
      \text{diameter}\; \{\o F(t_u,t_v,t_w)\,:\,
            u=i,i\!+\!1;\; v=j, j\!+\!1;\; w = k, k\!+\!1\! \} \;<\; 2/3.
\end{align}
for $0\leq i, j,k< N$. The nearness relations
(\ref{eq-nearby-third}) will make it possible to define a
continuous ternary operation $\o G$ that interpolates the $N^3$
discrete function values $\o F(t_i,t_j,t_k)$ ($0\leq i,j,k<N$).
Using the fact that these values obey (\ref{eq-majority}) we will
be able to make sure that the interpolated operation $\o G$ also
obeys (\ref{eq-majority}). Thus we will have $\la S^1,\o
G\ra\models \Sigma$ with $\o G$ continuous, in contradiction to
the known fact \cite{wtaylor-vohl} that $S^1\not\models\Sigma$;
this contradiction will complete the proof of the theorem.

Applying Lemma \ref{lem-23} to (\ref{eq-nearby-third}) we see that
for $0\leq i, j,k< N$ there is an arc $A_{ijk}$ of length $<2/3$
such that
\begin{align}               \label{eq-near-in-arc-third}
  \{\,  \o F(t_u,t_v,t_w) \,:\, u=i,i\!+\!1,\;\; v=j,j\!+\!1,\;\;
         w=k,k\!+\!1\, \} \,\;\subseteq\;  \,A_{ijk}.
\end{align}

Now the proof continues much like Part 2 of the proof of Theorem
\ref{th-idem-comm-23}; we omit the details. The function values
$\o F(t_i,t_j,t_k)$ will be interpolated to a continuous operation
$\o G$. The interpolation is done first along grid lines $\{ (t_i,
t_j, u)\}$, $\{(t_i,t,t_k)\}$ and $\{(s,t_j,t_k)\}$, where $s$,
$t$ and $u$ range over $S^1$. It is then extended to the grid
surfaces $\{ (t_i,t, u)\}$, $\{ (s, t_j, u)\}$ and
$\{(s,t,t_k)\}$, and finally to the entire $3$-dimensional figure
$(S^1)^3$. As before, it is carried out one cell at a time in the
given subdivision, and as before (\ref{eq-near-in-arc-third})
ensures that a continuous extension always exists, one cell at a
time.

To accommodate Equations  (\ref{eq-majority}) we need first
notice, for example, that $\o G(t_i,t_j,t_j) = \o
F(t_i,t_j,t_j)=t_j$, and likewise $\o G(t_i,t_{j+1},t_{j+1}) = \o
F(t_i,t_{j+1},t_{j+1})=t_{j+1}$. Therefore for $t$ ranging over
the arc $\o{t_j\,t_{j+1}}$ it is possible to define $\o G$ in such
a way that $\o G(t_i,t,t)=t$, which is a start on proving
(\ref{eq-majority}) for $\o G$. This can then be incorporated into
the determination of the two-dimensional interpolation $\o
G(t_i,t,w)$, by interpolating over two triangles, as we did in the
proof of Theorem \ref{th-idem-comm-23}. At the three-dimensional
level we must divide a cube into two triangular prisms. We omit
further details.
\end{Sketch}

\subsubsection{Comment on the proofs of Theorems
                           \ref{thm-23}--\ref{th-s1-majority}.}
                            \label{sub-proofs-comment}
Theorems  \ref{thm-23}, \ref{th-idem-comm-23} and
\ref{th-s1-majority}, in \S\ref{sub-mult-zero-one},
\S\ref{sub-idem-comm} and \S\ref{sub-majority}, each conclude that
$\mu(S^1,\Sigma)=2/3$ for a certain theory $\Sigma$. The proofs
for $\mu(S^1,\Sigma)\geq2/3$ are essentially identical: each
involves interpolating a discontinuous operation over a fine grid,
and producing a continuous operation. (The proof of Theorem
\ref{thm-23} is not directly phrased this way, but it could easily
be rewritten to this form.) We are confident that this method
would extend to many more simple non-Abelian theories $\Sigma$,
perhaps all of them. (Perhaps one would need to invoke
\cite{wtaylor-eri} to satisfy $\Sigma$ continuously at the
cellular level.)

This common argument relies essentially on Lemma \ref{lem-23},
which allows each cell to be mapped into an interval, which is
topologically very feasible. We believe it will be possible to
find analogs to  Lemma \ref{lem-23} for higher dimensions (e.g.\
for $S^n$); in that case the method may extend to the study of
non-Abelian simple theories on $S^n$.

On the other hand the three proofs for $\mu(S^1,\Sigma)\leq2/3$
seem to have arisen ad hoc, on a completely case-by-case basis. To
remind the reader: each of these proofs involved the construction
of an algebra $(S^1,\o F)$ satisfying $\Sigma$ and with $\chi(\o
F)\leq2/3$. At this time there seems to be little scope for
extension of these methods to another set $\Sigma$ of equations,
or to spaces of higher dimension.

\subsection{{$\Sigma=$ commutative idempotent binary;
$A=S^2$}.}    \label{sub-s2} Here we begin to explore whether the
method of Theorems \S\S\ref{thm-23}--\ref{th-s1-majority} will
extend to other spaces. We first note that the two-dimensional
sphere $S^2$ is incompatible with the spaces that appear in those
theorems (\cite{wtaylor-sae}; see also \cite[\S3.2 and
\S3.2.1]{wtaylor-asi}). In fact we will sketch a proof of

\begin{theorem}               \label{thm-mu-s2-majority}
If\/ $\Sigma$ is the theory either of a binary operation with zero
and one (\S\ref{sub-mult-zero-one}), or of a symmetric idempotent
operation (\S\ref{sub-idem-comm}), or of a ternary majority
operation (\S\ref{sub-majority}), then $\mu(S^2,\Sigma)\geq2/3$.
\end{theorem}

Before sketching the proof, we will state an analog of Lemma
\ref{lem-23}. In our previous applications of Lemma \ref{lem-23},
the essential part of the conclusion is that $F$ lies in some
convex subset of $S^1$, i.e. an arc. Let us suppose that $S^2$ is
given the great-circle metric, with diameter scaled to $1$. For a
subset $A\subseteq S^2$, we say that $A$ is {\em convex} iff for
each two points $P,Q\in A$, we have $d(P,Q)<1$ and
$\o{P\,Q}\subseteq A$.

\begin{lemma}        \label{lem-23-convex}
If\/ $F$ is a finite subset of $S^2$ with\/
$\text{diameter}(F)<2/3$, then there is a convex subset\/ $A$ of
$S^2$ such that\/ $F\subseteq A$.\ENDPROOF
\end{lemma}

The proof of Lemma \ref{lem-23-convex} is like that of Lemma
\ref{lem-23}, and omitted for now. Notice again that $2/3$ is best
possible for this conclusion: three equally spaced points on a
great circle form a set F of diameter $2/3$ that does not lie in a
convex set.

We will use two facts about convex subsets: the first is that the
intersection of any family of convex subsets is convex. The second
is the property of being an {\em absolute extensor} ({\bf AE}). A
metrizable space $A$ is defined to be an AE in the family of
metrizable spaces iff it satisfies the following property: if $B$
is a closed subspace of a metrizable space $F$, and if $g\FROM B
\TO A$ is a continuous function, then there exists a continuous
function $\phi\FROM F\TO A$ such that $\phi\!\upharpoonright\!
B=g$. (See e.g. \cite[pages 34--35]{hu}.) Each convex subset of
$S^2$ is homeomorphic to a convex subset (in the ordinary sense)
of the plane, and hence is an AE, by \cite[pages 84--87]{hu}. (See
also \cite{borsuk-retracts} for the general theory of AE's (and
absolute retracts).)\vspace{0.1in}

\hspace{-\parindent}%
{\em Sketch of proof of Theorem \ref{thm-mu-s2-majority}.} We will
restrict our attention to the case where $\Sigma$ is the theory of
a symmetric idempotent operation (\S\ref{sub-idem-comm}). For a
proof by contradiction, we assume that $\mu(S^1,\Sigma)\;<\;2/3$.
By (\ref{eq-def-mu}), there is an algebra $\mathbf A=\la S^2,\o
F\ra$ such that $\chi(\o F)<2/3$ and such that $\mathbf
A\models\Sigma$. As in the proofs of Theorems \ref{thm-23} and
\ref{th-idem-comm-23}, we give $(S^2)^2$ the sum metric. As
before, there exists $\delta>0$ such that $d((a,b),(c,d))<\delta$
implies $d(\o F(a,b),\o F(c,d))<2/3$.

Now let us assume that $(S^2)^2$ has been triangulated in such a
way that each $4$-simplex has diameter $< \delta$. Moreover the
triangulation must be symmetric in the following sense. Let
$\iota$ be the involution of $(S^2)^2$ given by
$\iota(a,b)=(b,a)$, where $a,b\in S^2$. Our symmetry condition is
that if $\sigma$ is a simplex of the triangulation, then so is
$\iota[\sigma]$. Our final condition is that the diagonal of
$(S^2)^2$ --- namely $\{(a,a)\,:\,a\in S^2\}$ --- must be a
subcomplex of this triangulation. Such a triangulation is clearly
possible.

We now proceed to define a continuous binary operation $\o G$ on
$S^2$, which is symmetric and idempotent. This will contradict the
known fact \cite[Theorem 1]{wtaylor-sae} that no such $\o G$
exists; this contradiction will complete the proof of Theorem
\ref{thm-mu-s2-majority}.

We define $\o G$ on simplices of successively higher dimension.
For a 0-simplex (point) $P$ we simply define $\o G(P)=\o F(P)$;
then obviously $\o G$ is symmetric and idempotent at the level of
$0$-simplices.

For each $4$-simplex $\sigma$, Lemma \ref{lem-23-convex} yields a
convex subset $A_{\sigma}$ of $S^2$ such that $\o
F[\o{\sigma}]\subseteq A_{\sigma}$. (Here $\o{\sigma}$ denotes the
closure of $\sigma$, which is $\sigma$ together with all its
subsimplices.) From the symmetry of $\o F$, we may further take
the sets $A_{\sigma}$ so that $A_{\iota(\sigma)}= A_{\sigma}$ for
all $\sigma$. We will use the $A_{\sigma}$'s in defining $\o G$
over simplices of dimensions 1, 2, 3 and 4. For a simplex $\tau$
of any dimension $\leq 4$, we define
\begin{align}
           A_{\tau} \EQ \bigcap\, \{ A_{\sigma} :
           \dim(\sigma)=4;\,
                                        \o{\tau}\subseteq
                                        \o{\sigma}\}.
\end{align}
It is not hard to check that $A_{\tau}$ is a nonempty convex
subset of $S^2$ such that
\begin{align}                \label{eq-rho-tau-behavior}
                  \o F[\o{\rho}]\subseteq A_{\rho};\;\;\;
          \text{if $\rho\leq\tau$, then $A_{\rho}\subseteq
          A_{\tau}$.}
\end{align}

We will now show inductively that for $n=1,2,3,4$, it is possible
to define $\o G$ on the $n$-skeleton of our triangulation in such
a way that, $\o G[\o{\tau}]\subseteq A_{\tau}$ for each
$n$-simplex $\tau$. We prove this for $n=3$; the other cases are
similar. If $\tau$ is a 3-simplex, then $\o G$ has already been
defined on all $2$-simplices in the boundary of $\tau$. For each
boundary $2$-simplex $\rho$, we have $\o F[\o{\rho}]\subseteq
A_{\rho}$ by (\ref{eq-rho-tau-behavior}). Also by
(\ref{eq-rho-tau-behavior}), we know that each of the sets
$A_{\rho}$ is a subset of $A_{\tau}$. Therefore the
two-dimensional extension of $\o G$ maps the boundary of $\tau$
into $A_{\tau}$. Since  $A_{\tau}$ is convex, and hence an AE,
there is a continuous extension of $\o G$ from the closed
3-simplex $\o{\tau}$ into $A_{\tau}$.

If we consider two closed 3-simplices, $\o{\tau}$ and $\o{\tau'}$,
then their overlap consists of closed 2-simplices; hence the
extensions to  $\o{\tau}$ and $\o{\tau'}$ agree on this overlap,
Thus the union of all such extensions is a well-defined continuous
function as desired. The desired condition $\o
G[\o{\tau}]\subseteq A_{\tau}$ was automatically fulfilled as we
went along.
Continuing in this manner to 4-simplices, we obtain a continuous
operation $\o G$ that agrees with $\o F$ on all vertices of the
triangulation.

 It remains to see that this operation can be made
to satisfy idempotence and symmetry. As for idempotence, since $\o
F$ satisfies $F(x,x) \wavy x$, we can easily define $\o G(x,x)$ to
be $x$. This may be taken as the definition of $\o
G\upharpoonright\sigma$, for each simplex $\sigma$ of the diagonal
subcomplex. We already have that $\o F[\o{\sigma}]\subseteq
A_{\sigma}$ for all $\sigma$, including those on the diagonal.
Since $\o G=\o F$ on the diagonal, we also have the required
condition that $\o G[\o{\sigma}]\subseteq A_{\sigma}$.
Incorporating this special case into our definition of $\o G$, we
now have a continuous idempotent operation.

As for symmetry, we merely need, for each simplex $\sigma$, to
define $\o G$ on $\sigma$ and $\iota[\sigma]$ at the same time.
(If $\sigma$ is a simplex of the diagonal subcomplex, then
$\sigma=\iota[\sigma]$, and so this condition has already been
met.) Inductively, we may assume that $\o G = \o G\Compos\iota$ on
all the boundary simplices of $\sigma$. Thus we simply define $\o
G$ on $\sigma$ as we did above, and on $\iota[\sigma]$ we define
$\o G$ by the formula $\o G = \o G\Compos\iota$. Clearly all the
conditions are met, and we now have a continuous, symmetric,
idempotent binary operation on $S^2$. This contradiction completes
the proof of the theorem.\ENDPROOF

\subsubsection{Comment on the proof of Theorem
         \ref{thm-mu-s2-majority}.}

In some ways the proof of Theorem \ref{thm-mu-s2-majority} may be
more comprehensible than those that we have supplied for Theorems
\ref{thm-23}, \ref{th-idem-comm-23} and \ref{th-s1-majority}, in
\S\ref{sub-mult-zero-one}, \S\ref{sub-idem-comm} and
\S\ref{sub-majority}. In those proofs we supplied a grid, which is
tantamount to a triangulation, but we needed to work with details
of that grid (often speaking, for instance, of $i$ and $i+1$,
etc.). In our proof of Theorem \ref{thm-mu-s2-majority}, we use
the general and inclusive notion of triangulation, which can be
discussed without reference to the detailed configuration of a
given triangulation.

It now seems right to conjecture that the method will go a lot
further than we have seen it here so far.

\subsection{An auxiliary theory.}
In 1986---see \cite[\S3.18, page 35]{wtaylor-cots}---we introduced
the following equational theory, known here as $\Sigma_1$:
\begin{align}                \label{eq-old-fix-p-a}
                 F(\phi^k(x),x,y)&\Wavy x \\
                     F(x,x,y)&\Wavy y, \label{eq-old-fix-p-b}
\end{align}
for $k\in\omega$, $k\geq 1$. We proved [{\em loc.~cit.}] that it
is incompatible with every compact Hausdorff space $A$ of more
than one element. In \cite[\S3.3.9]{wtaylor-asi} we proved that
$\Sigma_1$ has a $\lambda$-value
(\S\ref{sub-metric-approx-compat}) at least as large as
$\text{diameter}(A)/4$. Here we prove

\begin{theorem}
If $A$ is compact, then
   $\mu(A,\Sigma_1)\;\geq\; \text{diameter}(A)/2$\,.
\end{theorem}\begin{Proof}
To prove the theorem by contradiction, we may suppose that
$\mu(A,\Sigma_1)\,<\,\text{diameter}(A)/2$. In a manner by now
familiar, there exist (discontinuous) operations $\o F$ and $\o
{\phi}$ modeling (\ref{eq-old-fix-p-a}--\ref{eq-old-fix-p-b}) on
$A$, and positive real numbers
$\delta_0\leq\delta_1<\text{diameter}(A)/2$, such that $\o F$ and
$\o {\phi}$ are each constrained by $(\delta_0,\delta_1)$.

Since $A$ is compact, there exist $a,b\in A$ with
$d(a,b)\,=\,\text{diameter}(A)$. Choose arbitrary $q\in A$. By
compactness, the sequence $\o{\phi}{\,}^n(q)$ has a convergent
subsequence:
\begin{align*}
          \lim_{i\TO\infty}\o{\phi}{\,}^{n(i)}(q)\EQ c\in A .
\end{align*}
By the triangle inequality, either $d(b,c)\geq \text{diameter}
(A)/2$ or $d(a,c)\allowbreak\geq \allowbreak\text{diameter}
\allowbreak(A)/2$. Without loss of generality, we will assume that
 $d(b,c)\geq \text{diameter} (A)/2$. By \eqref{eq-old-fix-p-b},
\begin{align}         \label{eq-trick-seq}
   \o F\bigr(\o{\phi}{\,}^{n(i+1)}(q),\,\o{\phi}{\,}^{n(i)}(q),\,
                   b\bigr) \EQ \o{\phi}{\,}^{n(i)}(q)
\end{align}
for all $i$, and hence this sequence has $c$ as limit. On the
other hand, according to Lemma \ref{lem-limit-approx}, the
sequence in \eqref{eq-trick-seq} is eventually within $\delta_1$
of $\o F(c,c,b)\EQ b$ (by \eqref{eq-old-fix-p-a}). Therefore,
$d(b,c)\,\leq\,\delta_1\,<\, \text{diameter}(A)/2$, contrary to
our assumption. This contradiction completes the proof of the
theorem.
\end{Proof}

We note that in the proof the $(\delta_0,\delta_1)$-constraint on
$\o {\phi}$ was never used.

\subsection{A second auxiliary theory.}  \label{sub-aux-2}
In \cite[\S3.3.9]{wtaylor-asi} we introduced the following theory,
known here as $\Sigma_2$:
\begin{gather}         \label{eq-old-fix-p-a-bis}
           G(\psi_{m+k}(x,y),\,\psi_m(x,y),\,x,\,y) \Wavy x\\
             K(x,y)\Wavy G(u,u,x,y)\Wavy K(y,x),
                   \label{eq-old-fix-p-b-bis}
\end{gather}
for $m,k\in\omega$, with $k\geq1$. We proved [loc.~cit.] that
$\Sigma_2$ is incompatible with any compact $A$ with more than one
element. More precisely, we proved that $
\lambda_A(\Sigma_2)\leq\text{diameter}(A)/4$. Here we prove
something similar for $\mu$.

\begin{theorem}       \label{th-aux-for-log}
If $A$ is compact, then
   $\mu(A,\Sigma_2)\;\geq\; \text{diameter}(A)/2$\,.
\end{theorem}\begin{Proof}
To prove the theorem by contradiction, we may suppose that
$\mu(A,\Sigma_2)\,<\,\text{diameter}(A)/2$. In a manner by now
familiar, there exist (discontinuous) operations $\o G$, $\o K$
and $\o {\psi}_m$ modeling
(\ref{eq-old-fix-p-a-bis}--\ref{eq-old-fix-p-b-bis}) on $A$, and
positive real numbers $\delta_0\leq\delta_1<\text{diameter}(A)/2$,
such that $\o G$, $\o K$ and $\o {\psi}_m$ are each constrained by
$(\delta_0,\delta_1)$.

Let $a$ and $b$ be points of $A$ with $d(a,b)$ equal to the
diameter of $A$. By the triangle inequality, we have either
$d(a,\o K(a,b) \geq \text{diameter}(A)/2$ or $d(b,\o K(a,b) \geq
\text{diameter}(A)/2$. Without loss of generality, we shall assume
that
\begin{align}                      \label{eq-distance-b-K}
             d(b,\o K(a,b)) \;\geq\; \text{diameter}(A)/2.
\end{align}

Consider the sequence $\o\psi_i(b,a)$; by compactness it has a
convergent subsequence:
\begin{align*}
       \lim_{i\longrightarrow\infty} \o\psi_{n(i)}(b,a)\EQ c\in A.
\end{align*}
By \eqref{eq-old-fix-p-a-bis},
\begin{align}          \label{eq-trick-bis}
            \o G(\o\psi_{n(i+1)}(b,a),\,\o\psi_{n(i)}(b,a),\,b,\,a)
                        \EQ b
\end{align}
for all $i$. On the other hand, according to Lemma
\ref{lem-limit-approx}, the sequence in \eqref{eq-trick-bis} is
eventually within $\delta_1$ of $\o G(c,c,b,a)\EQ \o K(a,b)$ (by
\eqref{eq-old-fix-p-a-bis}). Therefore, $d(b,\o
K(a,b))\,\leq\,\delta_1\,<\, \text{diameter}(A)/2$, contrary to
\eqref{eq-distance-b-K}. This contradiction completes the proof of
the theorem.
\end{Proof}

Notice that the proof of Theorem \ref{th-aux-for-log} does not
mention the $(\delta_0,\delta_1)$-constraint on $\psi_m$, for any
$m$. Ignoring this constraint, we obtain the following sharper
version:

\begin{theorem}          \label{th-aux-for-log-better}
If\/ $A$ is a compact metric space of more than one element, then
there is no algebra\/ $\mathbf A\,=\,\la A,\o G, \o K,
\o{\psi}_m\ra_{m\in\omega}$ such that\/ $\chi(\o G)$ and\/
$\chi(\o K)$ are both\/ $\,<\,\text{\em diameter}(A)/2$, and\/
$\mathbf A\models\Sigma_2$.
\end{theorem}

\section{Dealing with composite operations.}
There may be a problem in carrying some of the results to
equations $\Sigma$ that involve composite operations. Suppose, for
example that $\o f$ is unary and $\chi(\o f)=\varepsilon$. If $s$
lies between $\o f(a)$ and $\o f(c)$ on a segment, then there
exists $b$ such that $\o f(b)$ lies with $\varepsilon$ of $s$. Our
equation of interest may, however, involve $\o g(\o f(b))$, and we
might like to know that this value is near to $\o g(s)$. With what
we have so far, we cannot conclude anything about the distance
between these two $\o g$-values.

\subsection{$n$-iterated jumps}
Let $(A,d)$ be a metric space, and $\mathbf A =(A,\o{F}_t)_{t\in
T}$ an algebra based on $A$. Recalling $\chi$ from
\S\ref{sub-cont-fail}, we define
\begin{align*}
        \chi_n(\mathbf A,d) \EQ \sup_{\tau}\; \chi(\o {\tau}) .
\end{align*}
Where $\tau$ ranges over all terms in operation symbols $F_t$ that
have depth $\leq n$, and where, for each $\tau$, $\o {\tau}$
denotes the term operation corresponding to $\tau$ in the algebra
$\mathbf A =(A,\o{F}_t)_{t\in T}$. We may also write
$\chi_{\infty}(\mathbf A,d)$ for the same supremum, taken over all
terms $\tau$.

 When the metric $d$ is clear
from the context, we may write $\chi_n(\mathbf A)$ for
$\chi_n(\mathbf A,d)$.

Finally, for $(A,d)$ a metric space, and $\Sigma$ a set of
equations of similarity type $\la n_t\,:\, t\in T\ra$, we define
\begin{align}                    \label{eq-def-mu-n}
         \mu_n(A,d,\Sigma) \EQ \inf\,\{\, \chi_n(\mathbf A,d) \,:\,
                 \mathbf A = (A,\o{F}_t)_{t\in T}
                    \models\Sigma\;\};
\end{align}
in other words, it is the infimum taken over all algebras built on
$A$ that satisfy $\Sigma$.  When the metric $d$ is clear from the
context, we may write $\mu_n(A,\Sigma)$ for $\mu_n(A,d,\Sigma)$.
We may also write $\mu_{\infty}(A,d,\Sigma)$ for the corresponding
infimum of $\chi_{\infty}$-values.

Obviously there is a uniform version
\begin{align*}
        \chi_{n}^u(\mathbf A,d) \EQ \sup_{\tau}\; \chi_u(\o {\tau})
        ,
\end{align*}
and likewise for $\mu_{n}^u$. Most of this paper deals with
compact metric spaces, on which the two concepts coincide, so we
will rarely mention $\chi_{n}^u$.

It is not hard to see that $\mu_n\leq\mu_{n+1}\leq\mu_{\infty}$
for all $n$, and moreover we generally expect that
$\mu_n<\mu_{n+1}<\mu_{\infty}$. Therefore, concerning estimates
from below, viz.\ $\mu_n>\varepsilon$, one should assert this for
$n$ as small as possible, in order to convey the most information.
On the other hand, such an estimate for a larger value of $n$ may
be all that is available, hence very valuable in itself.

\subsection{Iterated $(\delta,\varepsilon)$-closeness.}
                   \label{sub-iterated-closeness}
Let us say that a function $f\FROM A\TO B$ is {\em constrained by
$(\delta,\varepsilon)$}, or {\em
$(\delta,\varepsilon)$-constrained} iff it satisfies
\begin{align*}
       \text{if $d(x,x')<\delta$, then
       $d(f(x),f(x'))<\varepsilon$.}
\end{align*}
The notion is of course familiar, in that $f$ is defined to be
uniformly continuous iff for every $\varepsilon>0$ there exists
$\delta>0$ such that $f$ is $(\delta,\varepsilon)$-constrained.

In working with a finite direct power $A^n$ of a metric space
$(A,d)$, let us agree to give $A^n$ the following adjusted version
of the sum metric:
\begin{align}
        d\bigl((a_1,\ldots,a_n),\,(b_1,\ldots,b_n)\bigr) \EQ
              \frac{1}{n}\,\bigl(d(a_1,b_1)\,+\,\dots\,+\,d(a_n,b_n)\bigr).
\end{align}
This definition has the advantage that if $\text{diameter}(A)=1$,
then $\text{diameter}(A^n)=1$. It also figures in the detailed
proof of Lemma \ref{lem-constrained-chi} just below.

Now suppose that there are positive reals
$\delta_0,\delta_1,\ldots,\delta_n$ such that every operation of
$(A,F_t)_{t\in T}$ is constrained by the $n$ pairs
$(\delta_0,\delta_1),\,\allowbreak(\delta_1,\delta_2),\ldots,
\,\allowbreak(\delta_{n-1},\delta_n)$. In this case, we say that
{\em $(A,F_t)_{t\in T}$ is\/ $n$-constrained by\/}
$(\delta_0,\delta_n)$. The first lemma says that we may always
assume that the $\delta$'s form an increasing sequence.

\begin{lemma}      \label{lem-constrained-downward}
If\/ $\mathbf A=(A,F_t)_{t\in T}$ is\/ $n$-constrained by
$(\delta_0,\delta_n)$, then there are positive reals $\delta_i'$
(for $0\leq i\leq n$) such that $\delta_n'=\delta_n$, such that
$\delta_n'\geq \delta_{n-1}'\geq\cdots\geq\delta_1'\geq\delta_0'$,
and such that $(A,F_t)_{t\in T}$ is constrained by the $n$ pairs
$(\delta_0',\delta_1'),\,\allowbreak(\delta_1',\delta_2'),\ldots,
\,\allowbreak(\delta_{n-1}',\delta_n')$.
\end{lemma}\begin{Proof}
If the given $\delta_i$ do not already form a monotone increasing
sequence, then for some $i$ we have
\begin{align*}
         \delta_i > \delta_{i+1}\leq \delta_{i+2}\leq\cdots\leq
                     \delta_{n-1}\leq\delta_n.
\end{align*}
Let us define
\begin{align*}
    \delta_0'\EQ\delta_1'\EQ \cdots\EQ \delta_i'\EQ \delta_{i+1}\\
            \delta_j'\EQ \delta_j\quad(\text{for $i<j\leq n$}).
\end{align*}
It is clear that these values of $\delta_j'$ have the required
properties.
\end{Proof}

\begin{lemma}      \label{lem-constrained-chi}
If\/ $\mathbf A=(A,F_t)_{t\in T}$ is\/ $n$-constrained by
$(\delta_0,\delta_n)$, then $\chi_n(A,F_t)_{t\in
T}\leq\delta_n$.\ENDPROOF
\end{lemma}

We then define
\begin{align}
      \chi_n^{\star}(\mathbf A)&\EQ \inf\, \{\delta_n:
             (\exists \delta_0>0)\; \text{$\mathbf A=(A,F_t)_{t\in T}$
             is\/ $n$-constrained by $(\delta_0,\delta_n)$}\,\};
                     \label{eq-def-chi-n-star}\\
      \mu_n^{\star}(A,\Sigma) &\EQ \inf\,\{\, \chi_n^{\star}
           (\mathbf A,d) \,:\, \mathbf A = (A,\o{F}_t)_{t\in T}
                    \models\Sigma\;\}. \label{eq-def-mu-n-star}
\end{align}

Lemma \ref{lem-constrained-chi} then implies the first inequality
of

\begin{lemma}
   $\chi_n(\mathbf A)\,\leq\,\chi_n^{\star}(\mathbf A)$ and\/
  $\mu_n(\mathbf A,\Sigma) \,\leq\, \mu_n^{\star}(\mathbf
  A,\Sigma)$. If\/ $A\models\Sigma$, then these last two
  $\mu$-values are both zero.
\end{lemma}

In the sections that follow, we will be able to prove that
$\mu_n^{\star}(A,\Sigma)\geq K$ for certain $A$, $\Sigma$ and
$K>0.$ While this information is obviously less informative than
it would be to have $\mu_n(A,\Sigma)\geq K$, it nevertheless has
the virtues of being provable and of being a non-trivial
quantitative version of $A\not\models\Sigma$. In one case (see
\S\ref{sub-revisit-inj-binary}) we have
$\mu_2^{\star}(A,\Sigma)\geq K$ while $\mu(A,\Sigma)=0$. In this
case, $\mu_2^{\star}$ obviously conveys the greater amount of
information.

\subsection{Some consequences of      \label{sub-close-conseq}
                          $(\delta,\varepsilon)$-closeness.}
\begin{lemma}         \label{lem-limit-approx}
{\em(Limit theorem, approximate version.)} $f(x_i)$ approaches
$f(\lim x_i)$ within $\varepsilon$.
\end{lemma}

\begin{lemma}
{\em(Intermediate Value Theorem, approximate version.)} (Move here
from Lemma \ref{lem-approx-ivt}, \S\ref{sub-Y}.)
\end{lemma}

\begin{lemma}        \label{th-brouwer}
{\em (Brouwer Fixed-Point Theorem, approximate version.)}
\end{lemma}

\begin{lemma}           \label{th-borsuk-ulam-approx}
{\em (Borsuk-Ulam Theorem, approximate version.)}
\end{lemma}

\begin{lemma}             \label{lem-approx-into-circle}
Suppose that $A$ is a triangulable compact metric space (i.e.\ the
geometric realization of a finite simplicial complex). Let\/ $S^1$
be the ordinary $1$-sphere with arc-length distance, scaled to
have diameter $1$. Suppose that\/ $\o F\FROM A\TO S^1$ is
$(\delta, \varepsilon)$-constrained, where $0<\delta$ and\/
$0<\varepsilon<2/3$. Then there exists a continuous function $\o
G\FROM A\TO S^1$ such that $d(\o F(a),\o G(a))\,<\,\varepsilon$
for all $a\in A$.
\end{lemma}\begin{Sketch} The proof is much like that of Theorems
\ref{thm-23}, \ref{th-idem-comm-23} and \ref{th-s1-majority}, in
\S\ref{sub-mult-zero-one}, \S\ref{sub-idem-comm} and
\S\ref{sub-majority}, and especially like that of Theorem
\ref{thm-mu-s2-majority} in \S\ref{sub-s2} (even though this last
result is officially about the $2$-sphere).
\end{Sketch}

\begin{corollary}      \label{cor-comp-lambda-mu-simple-circle}
Suppose that $\Sigma$ is a simple theory, and $\mu(S^1,\Sigma)
\,<\,2/3$. Then $\lambda_{S^1}(\Sigma)\,\leq\,2\,\mu(S^1,\Sigma)$.
\end{corollary}

\begin{corollary}        \label{cor-comp-lambda-mu-k-circle}
{\em (Conjectured.)} Suppose that each equation of $\Sigma$
equates two terms of depth no more that $k$, and that
$\mu_k^{\star}(S^1,\Sigma) \,<\,2/3$. Then
$\lambda_{S^1}(\Sigma)\,\leq\,2\,\mu_k^{\star}(S^1,\Sigma)$.
\end{corollary}

Results like Corollaries
\ref{cor-comp-lambda-mu-simple-circle}--\ref{cor-comp-lambda-mu-k-circle}
must be relatively abundant. We will extend their range to other
spaces as tools become available.

\subsection{Revisiting \S\ref{sub-aux-2}: lattice-ordered groups.}
             \label{sub-lgroups}
Following \cite[\S3.3.10]{wtaylor-asi} we define $\Lambda\Gamma$
to be the following (doubly infinite) set of equations:
\begin{align}
   x&\Wavy x\wedge[(z_{m+k} - z_m) \,+\, (x\wedge y)]
                                         \label{eq-lg-c}\\
        x\wedge y &\Wavy   x\wedge[(u\,-\,u) + (x\wedge y)]
                   \Wavy y\wedge x,
                                            \label{eq-lg-a}
\end{align}
where $z_n$ ($n\in\omega$) are terms defined recursively as
follows:
\begin{align*}        
        z_0 \EQ 0;\quad\quad\quad  z_{n+1}\EQ (z_n\,+\,(x-(x\wedge y))).
\end{align*}
In \cite[\S3.3.10]{wtaylor-asi} we gave an easy proof that {\em
lattice-ordered groups} satisfy (\ref{eq-lg-c}--\ref{eq-lg-a}); in
other words Equations (\ref{eq-lg-c}--\ref{eq-lg-a}) are among the
consequences of the equational axioms of lattice-ordered group
theory (which we do not state here in detail). Thus any result of
the form $\mu_n(A,\Lambda\Gamma)\geq K$ or
$\mu_n^{\star}(A,\Lambda\Gamma)\geq K$---such as Theorem
\ref{th-lgroups} just below---implies the same result for the
theory of lattice-ordered groups.

The incompatibility of compact Hausdorff spaces with
lattice-ordered groups was proved by M. Ja.\ Antonovski\u\i\ and
A. V. Mironov \cite{anton-mir} in 1967. For compact metric spaces,
a positive value for $\lambda_A(\Lambda\Gamma)$ was established by
W. Taylor [{\em loc.~cit.}\/]. Here in \S\ref{sub-lgroups} we
prove a positive value for $\mu_3(A,\Lambda\Gamma)$.

Our method for estimating $\mu_{3}(A,\Lambda\Gamma)$ is to connect
$\Lambda\Gamma$ with the $\Sigma_2$ appearing in Equations
(\ref{eq-old-fix-p-a-bis}--\ref{eq-old-fix-p-b-bis}) of
\S\ref{sub-aux-2}. Lemma \ref{lem-interp-log} below will establish
an {\em interpretation}\footnote%
{At some point it may become appropriate to add a section on the
persistence of $\mu$-values under interpretation.} %
(in the sense of \cite{neumann,ogwt-mem}) of $\Sigma_2$ in
$\Lambda\Gamma$. (Thus $\Sigma_2$ is {\em a fortiori\/}
interpretable in lattice-ordered groups.) For every algebra
$\mathbf A \EQ\la A, \mm, \jj, \boxplus, \boxminus\ra$ in the
similarity type of $\Lambda\Gamma$, we define a new algebra
$\mathbf A'=\la A, \o G, \o K,\o{\psi}_m\ra_{m\in\omega}$, in the
similarity type of $\Sigma_2$, as follows. For $a,b,c,d\in A$, we
let
\begin{align*}
 \o G(a,b,c,d)&\EQ c\mm[(a\boxminus b)\boxplus(c\mm d)],\\
 \o
K(a,b)&\EQ a\mm b\\
\o{\psi}_m(a,b)&\EQ \o{z_n}(a,b),
\end{align*}
where $z_n$ is as above, and $\o{z_n}$ is the term operation
associated to $z_n$.

As noted above, the following lemma and theorem hold {\em a
fortiori} for lattice-ordered groups.

\begin{lemma}        \label{lem-interp-log}
If\/ $\mathbf A$ satisfies\/ $\Lambda\Gamma$, then\/ $\mathbf A'$
satisfies $\Sigma_2$.
\end{lemma} \begin{Proof}
We need to see that Equations \eqref{eq-old-fix-p-a-bis} and
\eqref{eq-old-fix-p-b-bis} hold in $\mathbf A'$. We look at
\begin{align}
G(\psi_{m+k}(x,y),\,\psi_m(x,y),\,x,\,y) \Wavy x
             \tag{\ref{eq-old-fix-p-a-bis}}
\end{align}
in detail. To prove its satisfaction in $\mathbf A'$, we need to
substitute  our definitions of $\o G$ and $\o{\psi}_m$ into
\eqref{eq-old-fix-p-a-bis} and verify the resulting equation under
$\Lambda\Gamma$. The reader may check that the resulting equation
is tantamount to \eqref{eq-lg-c}, which is one of the defining
equations of $\Lambda\Gamma$. Thus \eqref{eq-old-fix-p-a-bis}
holds in $\mathbf A'$. The proof for \eqref{eq-old-fix-p-b-bis} is
similar.
\end{Proof}

\begin{theorem}        \label{th-lgroups}
If\/ $A$ is a compact metric space of more than one point, then
$\mu_3(A,\Lambda\Gamma)\,\geq\, \text{diameter}(A)/2$.
\end{theorem}\begin{Proof}
To prove the theorem by contradiction, we may suppose that
$\mu_3(A,\Lambda\Gamma)\,<\,\text{diameter}(A)/2$.
By the definition \eqref{eq-def-mu-n}, there exist
(discontinuous) operations $\mm$, $\jj$, $\boxplus$, $\boxminus$
on $A$ such that $\mathbf A \EQ\la A, \mm, \jj, \boxplus,
\boxminus\ra$ satisfies $\Lambda\Gamma$, and such that
$\chi_3(\mathbf A)<\text{diameter}(A)/2$. This means
that
\begin{align}          \label{eq-every-log-tau}
\chi(\o{\tau})\;<\;\text{diameter}(A)/2
\end{align}
for every term-operation $\o{\tau}$ of $\mathbf A$ having depth
$\leq\,3$.

Now the algebra $\mathbf A'$ of Lemma \ref{lem-interp-log} is a
model of $\Sigma_2$, whose operations $\o G$ and $\o K$, each
being a term-operation of depth $\leq\,3$, both have $\chi$-value
$<\,\text{diameter}(A)/2$. This contradiction to Theorem
\ref{th-aux-for-log-better} completes the proof of our theorem.
\end{Proof}

\subsection{Revisiting \S\ref{subsub-inf-dim}: the injective
                                       binary operation.}
           \label{sub-revisit-inj-binary}
\subsubsection{$A=[0,1]$.}
We return our attention to Equations (\ref{eq-squared}) of
\S\ref{subsub-inf-dim}, which we repeat here for convenience:
\begin{align}        \tag{\ref{eq-squared}}
                F_0(G(x_0,x_1))\WAVY x_0,\quad\quad
                F_1(G(x_0,x_1)) \WAVY x_1.
\end{align}
Moreover, we again let $A=[0,1]$ with the ordinary Euclidean
metric. In \S\ref{subsub-inf-dim} we proved that
$\mu(A,\Sigma)=0$. Here we shall prove that
$\mu_2^{\star}(A,\Sigma)=1.$ In fact, we shall prove it in a
somewhat broader context.

\begin{theorem}         \label{th-inj-bin-star2}
Let $A=[0,1]$ be given any metric that induces the usual topology.
Then $\mu_2^{\star}(A,\Sigma)\EQ\text{\em diameter}(A)$.
\end{theorem} \begin{Proof}
We note first that clearly
$\mu_2^{\star}(A,\Sigma)\;\leq\;\text{diameter}(A)$ for any $A$
and any $\Sigma$. Thus to prove the theorem by contradiction, we
may suppose that $\mu_2^{\star}(Y,\Sigma)\,<\,\text{diameter}(A)$.
By Definitions (\ref{eq-def-chi-n-star}--\ref{eq-def-mu-n-star})
there exist (discontinuous) operations $\o F_0$, $\o F_1$ and $\o
G$ modeling (\ref{eq-squared}) on $A$, and positive real numbers
$\delta_0\leq\delta_2<\text{diameter}(A)$ such that $(A,\o F_0,\o
F_1,\o G)$ is 2-constrained by $(\delta_0,\delta_2)$. Thus there
exists a further positive real $\delta_1$ such that
\begin{align}               \label{eq-const-inj-bin}
   \text{$\o F_0$, $\o F_1$ and $\o G$ are each constrained by
   $(\delta_0,\delta_1)$ and by $(\delta_1,\delta_2)$.}
\end{align}

Since $[0,1]$ is compact, there exist $a_0,a_1\in A$ with
$d(a_0,a_1)=\text{diam}(A)$. For flexibility of notation, we take
two such pairs: $d(a_0,a_1)=d(b_0,b_1)=\text{diam}(A)$.
Considering the four real numbers
\begin{align*}
     \o G(a_0,b_0),\;\;\o G(a_1,b_0),\;\;\o G(a_0,b_1),\;\;
                 \o G(a_1,b_1),
\end{align*}
we may assume, without loss of generality, that the smallest among
them is $\o G(a_0,b_0)$. Again without loss of generality, we may
assume that $\o G(a_1,b_0))\,\leq\,\o G(a_0,b_1).$ In other words,
we have
\begin{align*}
             \o G(a_0,b_0)\;\leq\; \o G(a_1,b_0)\;\leq\;\o
             G(a_0,b_1).
\end{align*}
Thus, along the segment $\overline {(a_0,b_0)(a_0,b_1)}$ in the
square $[0,1]^2$, the $(\delta_0,\delta_1)$-constrained function
$\o G$ takes values that are above and below the value $\o
G(a_1,b_0)$. By Lemma \ref{lem-approx-ivt}, there exists $e\in
[0,1]$ such that
\begin{align*}
                  d\,(\o G(a_0,e),\, \o G(a_1,b_0))\;<\; \delta_1.
\end{align*}
For the $(\delta_1,\delta_2)$-constrained function $\o F_0$ we now
calculate, using  $\Sigma$:
\begin{align*}
     d(a_0,a_1)\EQ d\bigl(\o F_0(\o G(a_0,e)),\,\o F_0(\o G(a_1,b_0))\bigr)
         \;<\; \delta_2\;<\;\text{diameter}(A).
\end{align*}
This contradiction to our choice of $a_0, a_1$ completes the
proof.
\end{Proof}

We notice that in this proof we needed the
$(\delta_0,\delta_1)$-constraint only for the binary operation $\o
G$, and the $(\delta_1,\delta_2)$-constraint only for the unary
operations $\o F_0$, $\o F_1$. (In other words,
(\ref{eq-const-inj-bin}) contains more information than
necessary.) It would thus be possible to give Theorem
\ref{th-inj-bin-star2} a slightly sharper statement by modifying
the hypotheses according to this observation. Similar remarks
apply elsewhere in the paper. As far as we can see for now, such
an endeavor merits neither the effort involved nor the cumbersome
statements that would result.

\subsubsection{Comments on the proof of Theorem \ref{th-inj-bin-star2}}
Our estimate is made for $\chi_2^{\star}$ only. This proof does
not yield information on $\chi_2$. The reason is that we must be
able to estimate the effect of applying $\o F_0$, $\o F_1$ and $\o
G$ to the number $e$ that is supplied by Lemma
\ref{lem-approx-ivt}. Such an $e$ is not necessarily\footnote%
{Objection: if we look at the proof of Lemma \ref{lem-approx-ivt},
we see that $e$ really is in the range. This needs to be sorted
out before publication.} %
 in the range of our
operations, so that we cannot make the necessary estimate simply
by applying some term-operation $\o{\tau}$.

Comparing this proof with the corresponding proof for $\lambda$
that appears in \cite{wtaylor-asi}, we note a lot of similarity.
In fact this proof is the same almost verbatim.

\subsubsection{$A=[0,1]^2$.}

Once again, we work with these equations:
\begin{align}        \tag{\ref{eq-squared}}
                F_0(G(x_0,x_1))\WAVY x_0,\quad\quad
                F_1(G(x_0,x_1)) \WAVY x_1.
\end{align}
We shall suppose that the usual topology of $[0,1]$ is given by a
metric $d_0$ with the property that $d_0(0,1)\geq1$. We then let
$A=[0,1]^2$ with the metric defined as a sum (taxi-metric,
$L_1$-norm): $d((a,b),(c,d))=d_0(a,c)+d_0(b,d)$.

\begin{theorem}         \label{th-inj-bin-star2}
$\mu_2^{\star}(A,\Sigma)\,\geq\,1.$
\end{theorem} \begin{Proof}
To prove the theorem by contradiction, we may suppose that
$\mu_2^{\star}(A,\Sigma)\,<\,1$. In a manner by now familiar,
there exist (discontinuous) operations $\o F_0$, $\o F_1$ and $\o
G$ modeling (\ref{eq-squared}) on $A$, and positive real numbers
$\delta_0\leq\delta_1\leq\delta_2<1$ such that
\begin{align}               \label{eq-const-inj-bin}
   \text{$\o F_0$, $\o F_1$ and $\o G$ are each constrained by
   $(\delta_0,\delta_1)$ and by $(\delta_1,\delta_2)$.}
\end{align}

Now $B^2\EQ[0,1]^4$, and so the boundary of this space is a
three-sphere $S^3$. Let us consider the action of $\o G$ on this
three-sphere. Since $\o G$ is $(\delta_0,\delta_1)$-constrained,
it takes on $\delta_1$-close values at two antipodal points, by
our version of the Borsuk-Ulam Theorem (Theorem
\ref{th-borsuk-ulam-approx}). Without loss of generality, two
antipodal points have the form $((0,x_1),(y_0,y_1))$ and
$((1,u_1),(v_0,v_1))$. We thus have
\begin{align*}
           d\bigl(\,\o G\bigl((0,x_1),(y_0,y_1)\bigr),\,
             \o G\bigl((1,u_1),(v_0,v_1)\bigr)\,\bigr) \;<\;
             \delta_1\,.
\end{align*}
Since $\o F_0$ is $(\delta_1,\delta_2)$-constrained, Equations
$\Sigma$ yield
\begin{align*}
       1 \;&\leq\; d\bigl((0,x_1),(1,u_1)\bigr)\\
          &\EQ  d\bigl(\,\o F_0\o G\bigl((0,x_1),(y_0,y_1)\bigr),\,
             \o F_0 \o G\bigl((1,u_1),(v_0,v_1)\bigr)\,\bigr)
                                \:<\;\delta_2.
\end{align*}
\end{Proof}

\subsection{Group theory on spaces with the fixed-point property.}
In this section we let $\Gamma$ stand for any equational theory
whose models are groups. (Some variation is possible in choice of
primitive operations and axioms, but any such theory will do.) We
will assume that binary $+$ and unary $-$ are available, either as
primitives or as derived operations.

A will be a metric space that has the {\em fixed-point property}:
if $f\FROM A\TO A$ is continuous, then there exists $e\in A$ such
that $f(e)=e$. Until we know the full scope of Theorem
\ref{th-brouwer}, we will state and prove it only for a power
$[0,1]^n$, which is to say, for an $n$-simplex. A corresponding
result for $\lambda$ was proved in \S3.3.1 of \cite{wtaylor-asi}.

\begin{theorem}         \label{th-simplex-gp-star2}
Let $A=[0,1]^n$ be given any metric that induces the usual
topology, and let $\Gamma$ denote group theory. Then
$\mu_2^{\star}(A,\Gamma)\EQ\text{\em diameter}(A)$.
\end{theorem} \begin{Proof}
We note first that clearly $\mu_2^{\star}(A,\Gamma)\leq
\text{diameter}(A)$ for any $A$ and any $\Gamma$. Thus to prove
the theorem by contradiction, we may suppose that
$\mu_2^{\star}(A,\Gamma)\,<\,\text{diameter}(A)$. By Definitions
(\ref{eq-def-chi-n-star}--\ref{eq-def-mu-n-star}) there exist
(discontinuous) group operations $\boxplus$ and $\boxminus$ on
$A$, and positive real numbers $\delta_0\leq\delta_2
<\text{diameter}(A)$, such that $(A,\boxplus,\boxminus)$ is
2-constrained by $(\delta_0,\delta_2)$. Thus there exists a
further positive real $\delta_1$ such that
\begin{align}               \label{eq-const-inj-bin}
   \text{$\boxplus$ and $-$ are each constrained by
   $(\delta_0,\delta_1)$ and by $(\delta_1,\delta_2)$.}
\end{align}

Since $A$ is compact, there are points $a_0$, $a_1\in A$ with
$d(a_0,a_1)=\text{diameter}(A)$. Consider the function $f\FROM
A\TO A$ defined by $f(x)=(a_1\boxminus a_0)\boxplus x$. Since $f$
is $(\delta_0,\delta_1)$-constrained, Theorem (\ref{th-brouwer})
yields $e\in A$ such that $d(e,f(e))<\delta_1$. Now let $g\FROM
A\TO A$ be defined by $g(x)=x \boxplus(\boxminus e\boxplus a_0)$.
Since $g$ is $(\delta_1,\delta_2)$-constrained, we have
\begin{align*}
           d(a_0,a_1)\EQ d\bigl(g(e),g(f(e))\bigr) \;<\;
                            \delta_2 \;<\; \text{diameter}(A).
\end{align*}
This contradiction to the choice of $a_0$ and $a_1$ completes the
proof of the theorem.
\end{Proof}

For unary operations of the form $x\GOESTO x+a$ and $x\GOESTO
x-b$, the constraints in \eqref{eq-const-inj-bin} are redundant,
since $x-b$ is the same as $x+a$, where $a=-b$. (For the full
binary operations, they may not be redundant.) This redundancy may
be seen in the proof, in the fact that we applied the constraints
\eqref{eq-const-inj-bin} only to operations of the form $x+a$.
Thus \eqref{eq-const-inj-bin} turns out to contain more
information than is necessary for the proof.

\subsection{Groups of exponent N on $\Reals$.}
In this section we let $\Gamma_N$ stand for any equational theory
whose models are (additively written) groups satisfying
$x+\cdots+x\wavy0$ (where $x$ appears $N$ times on the left of
this equation). For $N=2$, this theory was known
\cite[\S3.3.7]{wtaylor-asi} to be incompatible with $\Reals$.

We will use the fact that any function $f\FROM\Reals\TO\Reals$
that cycles a set of $N$ elements must have an approximate fixed
point, by a minor variation on Theorem \ref{lem-approx-ivt}.

\begin{theorem}
Let $\Reals$ be given any metric that induces the usual topology,
and let\/ $\Gamma_N$ denote group theory with exponent N. Then
$\mu_2^{\star}(\Reals,\Gamma_N)\geq\text{\em radius\/}(A)$.
\end{theorem} \begin{Proof}
To prove the theorem by contradiction, we may suppose that
$\mu_2^{\star}(A,\Gamma_N)\,<\,\text{radius}(A)$. By Definitions
(\ref{eq-def-chi-n-star}--\ref{eq-def-mu-n-star}) there are a
(discontinuous) exponent-N group operation $\boxplus$ on $A$, and
positive real numbers $\delta_0\leq\delta_2 <\text{radius}(A)$,
such that $(A,\boxplus)$ is 2-constrained by
$(\delta_0,\delta_2)$. Thus there exists a further positive real
$\delta_1$ such that
\begin{align*}               
   \text{$\boxplus$\; is constrained by
   $(\delta_0,\delta_1)$ and by $(\delta_1,\delta_2)$.}
\end{align*}

Let $\o0$ be the unit element of the group $(\Reals,\boxplus)$.
Since $\delta_2<\text{radius}(A)$, there exists $a\in\Reals$ such
that $d(a,\o 0)>\delta_2$. We consider the function
$f\FROM\Reals\TO \Reals$ given by $f(x)=x\boxplus a$. Clearly
$f(\o0)=a$, $f(a)=2a$, $f(2a)=3a$, \ldots and $f((N-1)a)=\o0$.
Since $f$ is $(\delta_0,\delta_1)$-constrained, (a variant on)
Theorem \ref{lem-approx-ivt} yields $e\in \Reals$ such that
$d(e,f(e))<\delta_1$. Let $g(x)=(N-1)e\boxplus x$. Since $g$ is
$(\delta_1,\delta_2)$-constrained, we have
\begin{align*}
     d(\o 0,a) \EQ d(g(e), g(f(e))) \;<\; \delta_2,
\end{align*}
in contradiction to our choice of $a$. This contradiction
completes the proof of the theorem.
\end{Proof}

\subsection{$A=Y$, the  triode; $\Sigma=$ lattice theory.}
                                       \label{sub-Y}

 Let $A,B,C,D$ be four non-collinear points in
the Euclidean plane, with $D$ in the interior of $\Tri{A}{B}{C}$.
Our space $Y$ is defined to be the union of the three (closed)
segments $AD$, $BD$ and $CD$, called {\em legs}, with the topology
inherited from the plane. In fact, in order to give $Y$ a definite
metric $d$, we will further require that $\Tri{A}{B}{C}$ be
equilateral with $D$ at its center, and that each leg have unit
length. We then let $d$ be the metric of the plane, as inherited
by $Y$.

For \S\ref{sub-Y} we let $\Sigma$ consist of axioms for lattice
theory (expressed in terms of $\wedge$ and $\vee$). It was proved
by A. D. Wallace in the mid-1950's (see \cite[{\em Alphabet
Theorem}, page xx]{wallace} for a statement of the result) that
the triode $Y$ is not compatible with $\Sigma$. Taking
$\mu^{\star}_3$ as defined in \S\ref{sub-iterated-closeness}, we
shall prove the sharper result that

\begin{theorem}     \label{th-triode}
$\mu_3^{\star}(Y,\Sigma)\geq\,0.5$.
\end{theorem}
Before proving Theorem \ref{th-triode} we state and prove one
Lemma. It is our discontinuous approximate replacement for the
Intermediate Value Theorem.

\begin{lemma} \label{lem-approx-ivt}
 Suppose that $f$ maps a convex subset of $\Reals$
into $\Reals$, and that $f$ is $(\delta,\varepsilon)$-constrained
for some $\delta,\varepsilon>0$. If $a<c$ and $s$ is between
$f(a)$ and $f(c)$, then there exists $b$ with $a\leq b\leq c$ and
with\/ $d(f(b),s)<\varepsilon/2$.
\end{lemma} \begin{Proof}
Consider a finite sequence of reals that begins with $a$ and ends
with $c$, and such that every step is smaller than $\delta$. The
corresponding function-values take steps smaller than
$\varepsilon$ while traversing the interval between $f(a)$ and
$f(c)$. Moreover $s$ must lie in one of these $f$-intervals
smaller than $\varepsilon$; hence the conclusion.
\end{Proof}

\hspace{-\parindent}%
{\em Proof of Theorem \ref{th-triode}.}\quad For a contradiction,
suppose that $\mu_3^{\star}(Y,\Sigma)\,<\,0.5$. By Definitions
(\ref{eq-def-chi-n-star}--\ref{eq-def-mu-n-star}) there exist
(discontinuous) lattice operations $\mm$ and $\jj$ on $Y$, and
positive real numbers $\delta_0$ and $\delta_3$, such that
$(Y,\mm,\jj)$ is 3-constrained by $(\delta_0,\delta_3)$, and
moreover such that $\delta_3<0.5$. Thus there exist further
positive reals $\delta_1, \delta_2$ such that
\begin{align}
   \text{$\mm$ and $\jj$ are each constrained by
   $(\delta_0,\delta_1)$, by $(\delta_1,\delta_2)$, and by
   $(\delta_2,\delta_3)$.}
\end{align}
By Lemma \ref{lem-constrained-downward} we may assume that
$\delta_0\leq\delta_1\leq\delta_2\leq\delta_3$. \vspace{0.1in}

\hspace{-\parindent}%
{\bf Part 1.} We shall prove  that either $A\mm D$ or $A\jj D$
lies in the leg $AD$ (and similarly for $B$ and $D$). If $A\mm D$
does not lie in $AD$, then we have $D$ between $A=A\mm A$ and
$A\mm D$. Consider the function of meeting with $A$, viz.
$X\GOESTO X\mm A$. Since $D$ is between two of its values, we may
apply Lemma \ref{lem-approx-ivt} to obtain $E\in AD$ with $d(A\mm
E,D)< \delta_1$. Now, joining with $A$, we have $d(A\jj(A\mm
E),A\jj D)<\delta_2$; by $\Sigma$ this may be simplified to
$d(A,A\jj D)<\delta_2$. In other words $A\jj D$ lies in $AD$ as
desired. \vspace{0.1in}

\hspace{-\parindent}%
{\bf Part 2.} We shall prove that either $A\mm D$ or $A\jj D$ lies
within $\delta_2$ of $A$ (and similarly with $A$ changed to $B$
and to $C$). By Part 1, the three points $A$, $A\mm D$ and $A\jj
D$ lie along a segment. Without loss of generality we have $A\jj
D$ between $A=A\mm A$ and $A\mm D$ on that segment. If we consider
the function of meeting with $A$ (as in Part 1), then  Lemma
\ref{lem-approx-ivt} again yields $E$ such that $d(A\mm E,\,A\jj
D)<\delta_1$. As in Part 1, joining with $A$ again yields
$d(A,A\jj D)<\delta_2$.\vspace{0.1in}

\hspace{-\parindent}%
{\bf Part 3.} From Part 2, we may assume, without loss of
generality, that
\begin{align}
 d(A,A\jj D)<\delta_2 \text{\quad and \quad} d(B,B\jj D)<\delta_2.
                  \label{eq-joins-near}
\end{align}
(The two vertices might be $A$ and $C$ or $B$ and $C$, and both
operations might be meets rather than joins, but surely two of the
three end-vertices must have the same pattern.)\vspace{0.1in}

\hspace{-\parindent}%
{\bf Part 4.} $A\jj B$ cannot lie in both of the disjoint sets
$[A,D)$ and $[B,D)$ (these are two of the legs, minus the endpoint
$D$). Without loss of generality we will assume that $A\jj B$ is
not in $[B,D)$. Therefore $D$ lies between $B=B\jj B$ and $A\jj
B$. By a familiar argument (this time involving joining with $B$)
we obtain $d(D,E\jj B)< \delta_1$ for some $E$. Now meeting with
$B$, we have $d(B\mm D,B)<\delta_2$.\vspace{0.1in}

\hspace{-\parindent}%
{\bf Part 5.} Taking the conclusion of Part 4, and joining with D,
yields $d(D,\allowbreak D\jj B)\,\allowbreak=\,d(D\jj(B\mm D),D\jj
B)<\delta_3$. Combining this with (\ref{eq-joins-near}), we have
\begin{align*}
       d(B,D)\;\leq\; d(B,D\jj B)\,+\, d(D\jj B,B)
       \;\leq\; \delta_3 + \delta_2 \;<\; 1.
\end{align*}
Here we have a contradiction to the fact that $d(B,D)=1$ (see the
start of \S\ref{sub-Y}), which completes the proof of Theorem
\ref{th-triode}.\ENDPROOF

\subsubsection{Comments on the proof of Theorem \ref{th-triode}}
Our estimate is made for $\chi_3^{\star}$ only. This proof does
not yield information on $\chi_3$. The reason is that we must be
able to estimate the effect of applying $\mm$ and $\jj$ to the
points called $E$ in Parts 1, 2 and 4. Such an $E$ is not
necessarily in the range of our operations, so that we cannot make
the necessary estimate simply by applying some term-operation
$\o{\tau}$.

Comparing this proof with the corresponding proof for $\lambda$
that appears in \cite{wtaylor-asi}, we note a lot of similarity.
It seems as though we could work out a theory for a composite
measure, including the possibility of limited jumps and of
approximate satisfaction. This will have to await a later date.

\subsection{A very special space.}

Section under construction.

For $\alpha$ any real number with $0<\alpha<1$, we define
\begin{align*}
    A_{\alpha}\EQ \{\,(x,y,z)\,:\; &x^2+y^2=1 \;\&\\
        &(-\alpha x\leq z\leq\alpha x\;\;\text{or}\;\;
           z=0)\, \}\;\subseteq\;\Reals^3.
\end{align*}
This space may easily be sketched as a subset of a cylinder in
$\Reals^3$. We give it the rectangular or taxicab metric in that
space: $d(\mathbf x,\mathbf y)=\sum|x_i-y_i|$. (Notice that the
spaces $A_{\alpha}$ are all homeomorphic one to another, but the
homeomorphisms are not isometries.)

Notice that for $(x,y,z)\in A_{\alpha}$ with $x<0$, the definition
yields $z=0$ as the only possible value for $z$. Thus $A_{\alpha}$
contains the circle
\begin{align*}
           C\EQ\{\,(x,y, 0)\,:\, x^2 + y^2 = 1\,\},
\end{align*}
and for negative $x$, these are the only points in $A_{\alpha}$.
For positive $x$, there are other points $(x,y,z)$. The farthest
of these from the circle $C$ are $(1,0,\alpha)$ and
$(1,0,-\alpha)$. Thus $\alpha$ is a measure of how far
$A_{\alpha}$ extends away from the circle $C$.

For future reference, we define a closed curve $f$  in
$A_{\alpha}$ (for $0\leq t\leq 2\pi$), as follows:
\begin{align*}
    f(t) &\EQ \begin{cases}
              \;\; (\cos t,\,\sin t,\,-\alpha\cos t)  &\text{if
                           $\cos t\geq0$} \\
               \;\;(\cos t,\,\sin t,\,0)
                &\text{if $\cos t\leq0$}.
          \end{cases}
\end{align*}
($f$ maps, so to speak, to the lower periphery of $A_{\alpha}$.)
Concerning the point $
        A \EQ (1,0,\alpha)\in A_{\alpha},
$ We proved in \cite[{\em loc.\ cit.}]{wtaylor-asi} that $A$ has
distance at least $2\alpha$ from every point $B=(\cos t,\cdots)$
in the image of $f$.

For this section, when we refer to a closed curve, we mean a
continuous map with domain $S^1$. We view our $f(t)$ as such a
closed curve, by representing $S^1$ as $\Reals/2\pi$, and relying
on the periodicity of the trigonometric functions. Finally when we
say that closed curves $g_0(t)$ and $g_1(t)$ in $A$ are homotopic,
we mean that there exists a map $G\FROM S^1\times [0,1]\TO A$ such
that $G(t,i) = g_i(t)$ for $t\in S^1$ and $i\in\{0,1\}$.

For $\mu_2^{\star}(A_{\alpha},\Gamma)\geq\alpha$, we restrict
$\Gamma$ to be a set of equations in the operations $+$ and $-$
that contains the three equations
\begin{gather}   \label{eq-ext-group-thy}
          (x + y) + (-y) \WAVY x \\
              x + 0 \WAVY x; \quad \quad 0+x\WAVY x.
                  \label{eq-ext-group-thy-b}
\end{gather}
We shall prove that for such a $\Gamma$,
\begin{align}             \label{eq-A-alpha-upper}
       \mu_2^{\star}(A_{\alpha},\Gamma)\;\geq\;\alpha.
\end{align}
The proof is by contradiction. To this end, we assume now that
$\mu_2^{\star}(A_{\alpha},\Gamma)\,<\,\alpha$.

By definition of $\mu_2^{\star}$, there exist positive reals
$\delta_0\leq\delta_1\leq\delta_2$ and (discontinuous) operations
$\boxplus,\boxminus$ obeying $\Gamma$, with both operations
constrained by $(\delta_0,\delta_1)$ and by $(\delta_1,\delta_2)$.
Let $0=A_0, A_1, A_2,\ldots,A_k=A$ be a sequence of members of
$A_{\alpha}$, with $d(A_i,A_{i+1})<\delta_0$ for each $i$. By the
$(\delta_0,\delta_1)$-constraint, we have
\begin{align*}
          d(A_i\boxplus\! f(t),\,A_{i+1}\boxplus\! f(t)) \;<\; \delta_1
\end{align*}
for each $i$ and each $t$. By (a version of) Lemma
\ref{lem-approx-into-circle}, for each $i$ there is a continuous
curve $\o g_i\FROM S^1\TO A_{\alpha}$
\begin{align}             \label{eq-g-near-f(t)plus}
         d(\o g_i(t), A_i\boxplus\! f(t))\;<\;\delta_1
\end{align}
for each $i$ and each $t$. Combining the last two inequalities, we
have
\begin{align*}
         d(\o g_i(t),\o g_{i+1}(t))\;<\;3\,\delta_1
\end{align*}
for each $i$ and each $t$. Assuming now that $3\delta_1$ is less
than the diameter of a circle in our model, we know that all the
functions $g_i$ are homotopic. In particular $g_k$ is homotopic to
$f(t)$ and hence maps onto $\{f(t):\pi/2\leq t\leq3\pi/2\}$. By
surjectivity and \eqref{eq-g-near-f(t)plus} there exists $t_0$
such that
\begin{align}            \label{eq-here-ist0}
           d(A \boxplus \!f(t_0), \,(-1,0,0))\;<\; \delta_1.
\end{align}

By reasoning similar to that for  \eqref{eq-here-ist0}, except
using the first equation of (\ref{eq-ext-group-thy-b}), we have
successively nearby points $0=B_0,B_1,\ldots,B_m=f(t_0)$ and maps
$h_i$, each close to the corresponding $f(t)\boxplus B_i$. In this
way, we arrive at the existence of $s_0$ with
\begin{align}            \label{eq-here-iss0}
           d(f(s_0) \boxplus\! f(t_0),\, (-1,0,0))\;<\;\delta_1.
\end{align}
%
Now from \eqref{eq-here-ist0} and \eqref{eq-here-iss0}, from the
group equation \eqref{eq-ext-group-thy}, and from the
$(\delta_1,\delta_2)$-constraint, we have
\begin{align*}
   d\bigl(A,(-1,&0,0)\boxminus\! f(t_0)\bigr)\\ &\EQ
         d \bigl((A\boxplus\! f(t_0))\boxminus \!f(t_0),(-1,0,0)\boxminus\!
         f(t_0)\bigr) \;<\; \delta_2;\\
    d\bigl(f(s_0),(-1,&0,0)\boxminus\! f(t_0)\bigr)\\ &\EQ
         d \bigl((f(s_0)\boxplus\! f(t_0))\boxminus \!f(t_0),(-1,0,0)\boxminus\!
         f(t_0)\bigr) \;<\; \delta_2.
\end{align*}
By the triangle inequality, $d(A,f(s_0))<2\delta_2<2\alpha$. This
contradicts our earlier assertion that every point in the image of
$f$ is at least $2\alpha$ from $A$, and thus the the result is
proved.

Also note that $A_{\alpha}$ is compatible with H-space theory.
(For the moment this is left to the reader.)

ONE FINAL PIECE would be to work out the diameter of $A_{\alpha}$.
Then check out the range of normalized values of $\lambda$. It
looks like we would still get a large range of $\lambda$-values.

\input{jour}

\vspace{\fill} \hspace*{-\parindent}%
\parbox[t]{3.0in}{ Walter Taylor \\ Mathematics Department\\ University
of Colorado\\
 Boulder, Colorado \ 80309--0395\\ USA\\ Email:
 {\tt walter.taylor@colorado.edu}}
\end{document}

%% file: jour.tex
\newcommand{\AEQ}{{\em Aequationes Mathematicae}}
\newcommand{\AAM}{{\em Advances in Applied Mathematics}}
\newcommand{\AMS}{{American Mathematical Society}}
\newcommand{\AU}{{\em Algebra Universalis}}
\newcommand{\ANNALS}{{\em Annals of Mathematics}}
\newcommand{\AML}{{\em Annals of Mathematical Logic}}
\newcommand{\ANNALEN}{{\em Mathematische Annalen}}
\newcommand{\BAMS}{{\em Bulletin of the \AMS}}
\newcommand{\BPAN}{{\em Bulletin de l'Aca\-d\'{e}\-mie Po\-lo\-naise
    des Sciences, S\'{e}rie des sciences math., astr. et phys.}}
\newcommand{\CAMB}{{\em Proceedings of the Cambridge Philosophical Society}}
\newcommand{\CANAD}{{\em Canadian Journal of Mathematics}}
\newcommand{\COLLOQ}{{\em Col\-lo\-qui\-um Ma\-the\-ma\-ti\-cum}}
\newcommand{\COLLOQUIA}{{\em Col\-lo\-qui\-a Ma\-the\-ma\-ti\-ca
    Soc\-i\-e\-ta\-tis Bolyai J\'{a}nos}}
\newcommand{\DM}{{\em Discrete Mathematics}}
\newcommand{\EM}{{\em l'En\-seigne\-ment Ma\-th\'{e}\-matique}}
\newcommand{\FM}{{\em Fun\-da\-men\-ta Ma\-the\-ma\-ti\-cae}}
\newcommand{\HOUS}{{\em Hous\-ton Journal of Mathematics}}
\newcommand{\INDAG}{{\em In\-da\-ga\-ti\-o\-nes Ma\-the\-ma\-ti\-cae}}
\newcommand{\JALG}{{\em Journal of Algebra}}
\newcommand{\JAMS}{{\em Journal of the Australian Mathematical Society}}
\newcommand{\JPAA}{{\em Journal of Pure and Applied Algebra}}
\newcommand{\JPAL}{{\em Journal of Pure and Applied Logic}}
\newcommand{\JSL}{{\em J. Symbolic Logic}}
\newcommand{\LNM}{{\em Lecture Notes in Mathematics}}
\newcommand{\MONTHLY}{{\em American Mathematical Monthly}}
\newcommand{\MUSSR}{{\em Mathematics of the {\sc USSR} --- Sbor\-nik}}
\newcommand{\NAMS}{{\em Notices of the \AMS}}
\newcommand{\NORSK}{{\em Norske Vid. Selssk. Skr. I, Mat. Nat.
                Kl. Chris\-ti\-a\-nia}}
\newcommand{\ORD}{{\em Order}}
\newcommand{\PACIFIC}{{\em Pacific Journal of Mathematics}}
\newcommand{\PAMS}{{\em Proceedings of the \AMS}}
\newcommand{\REM}{{\em Research and Exposition in Mathematics}}
\newcommand{\SF}{{\em Semigroup Forum}}
\newcommand{\SCAND}{{\em Ma\-the\-ma\-ti\-ca Scan\-di\-na\-vi\-ca}}
\newcommand{\SIAMJC}{{\em {\sc Siam} Journal of Computing}}
\newcommand{\SZEGED}{{\em Acta Sci\-en\-ti\-a\-rum Ma\-the\-ma\-ti\-ca\-rum}
                (Sze\-ged)}
\newcommand{\TAMS}{{\em Transactions of the \AMS}}
\newcommand{\TCS}{{\em Theoretical Computer Science}}

\newcommand{\ALVIN}{{\em Algebras, Lattices, Varieties}}
\newcommand{\Wads}{{Wads\-worth and Brooks-Cole Publishing Company,
                    Monterey, CA}}
\newcommand{\NorthHolland}{{North-Holland Publishing Company, Amsterdam}}
\newcommand{\Birk}{{Birkh\"{a}user}}
\newcommand{\Garcia}{Garc\'{\i}a}

%% file: jumps.bbl
\begin{thebibliography}{99}
\bibitem{adams}
    J. F. Adams, On the non-existence of elements of
    Hopf-invariant one, \ANNALEN\ (2) 72 (1960), 20--104.
\bibitem{anderson}
    R. M. Anderson, ``Almost'' implies ``near,'' \TAMS\ 296 (1986),
    229--237.
\bibitem{anton-mir}
M. Ja.\ Antonovski\u\i\ and A. V. Mironov,  On the theory of
      topological $l$-groups. (Russian. Uzbek summary) Dokl. Akad.
      Nauk UzSSR 1967, no. 6, 6--8. MR 46 \#5528.
\bibitem{band-hed}
   H.-J. Bandelt and J. Hedl\'{\i}kov\'a, Median algebras, \DM\
   45 (1983), 1--30.
\bibitem{borsuk}
K. Borsuk, Drei S\"atze \"uber die n-dimensionale euklidische
Sph\"are, Fund. Math. 20 (1933), 177--190.
\bibitem{borsuk-retracts}
  \Dash      {\em Theory of retracts},
   Monografie Matematyczne, Tom 44, Pa\'nstwowe Wydawnictwo Naukowe,
    Warsaw, 1967, 251 pages.\bibitem{bott}
   R. Bott,
         On symmetric products and the Steenrod squares, \ANNALS\
         57 (1953), 579--590.
\bibitem{caviness}
B. F. Caviness and J. R. Johnson, eds., {\em Quantifier
       Elimination and Cylindrical Algebraic Decomposition.
       Texts and Monographs in Symbolic Computation},
         New York: Springer-Verlag, 1998.
\bibitem{choe}
T. H. Choe, On Compact Topological Lattices of Finite Dimension,
{\em Transactions of the American Mathematical Society} 140
(1969), 223--237.
\bibitem{clifford}
A. H. Clifford, Connected ordered topological semigroups with
      idempotent endpoints, I. \TAMS\ 88 (1958), 80--98.
\bibitem{collins}
G. E. Collins, Quantifier Elimination for the Elementary Theory
       of Real Closed Fields by Cylindrical Algebraic Decomposition,
       {\em Lecture Notes in Computer Science} 33 (1975), 134-183.
\bibitem{dug-gran}
 J. Dugundji and A. Granas, {\em Fixed-point theory I}, Pa\'nstwowe
 Wy\-daw\-nic\-two Naukowe, Warsaw, 1982.
\bibitem{faucett}
W. M. Faucett, Compact semigroups irreducibly connected between
    two points; Topological semigroups and continua with
    cutpoints. \PAMS\ 6 (1955), 741--756.
\bibitem{ogwt-mem}
O. C. \Garcia\ and W. Taylor,
       {\em The lattice
    of interpretability types of varieties}, {\sc Memoirs of the American
    Mathematical Society}, Number 305, iii+125 pages.  MR 86e:08006a.
\bibitem{hu}
S. T. Hu,   {\em Theory of retracts}, Wayne State University
Press, Detroit, 1965.
\bibitem{james}
        I. M. James,
        Multiplication on spheres, I, II, \PAMS\ 13
        (1957), 192--196 and \TAMS\ 84 (1957), 545--558.
\bibitem{kapl-amerj}
I. Kaplansky, Topological Rings,  {\em American Journal of
Mathematics} 69 (1947), 153--183.

\bibitem{koch-wall-a}
 R. J. Koch and A. D. Wallace,
          Admissibility of semigroup structures on continua, AMS
               Transactions 88 (1958), 277--287.
\bibitem{koch-wall-b}
    \Dash
       Topological semilattices and their underlying spaces,
               Semigroup Forum 1 (1970), 209--223.
\bibitem{matousek}
J. Matou\v sek, {\em Using the Borsuk-Ulam Theorem, Lectures on
   Topological Methods in Combinatorics and Geometry},
   Springer-Verlag, Berlin, 2000.
\bibitem{mckenzie-cubes}
R. McKenzie, On spectra, and the negative solution of the
       decision problem for identities having a non-trivial
       finite model, \JSL\ 41 (1975), 186--196.
\bibitem{neumann}
W. D. Neumann, On Mal'cev conditions, \JAMS\ 17 (1974), 376--384.
\bibitem{pour-richards}
M. B. Pour-El and J. I. Richards,  {\em Computability in
      analysis
     and physics.} Perspectives in Mathematical Logic.
    Springer-Verlag, Berlin, 1989. xii+206 pp.
\bibitem{spanier}
     E. Spanier, {\em Algebraic Topology}, McGraw-Hill, New York, 1966.
\bibitem{steinlein}
  H. Steinlein, Borsuk's  theorem and its generalizations
  and applications: a survey, pages 166--235 in A. Granas, ed.,
  {\em M\'ethodes topologiques en analyse non lin\'eaire
  (S\'eminaire scientifique OTAN)}, Les Presses de l'Universit\'e
  de Montr\'eal, 1985.
\bibitem{tarski-1931}
A. Tarski, Sur les ensembles d\'efinissables de nombres r\'eels,
\FM\ 17 (1931), 210--239.

\bibitem{tarski-1948}
\Dash A Decision Method for Elementary Algebra and Geometry, RAND
Corp. monograph, 1948.

\bibitem{tarski-1951}
  \Dash {\em A Decision Method for Elementary Algebra and Geometry},
2nd ed. Berkeley, CA: University of California Press, 1951.

\bibitem{wtaylor-cmc}
W. Taylor,
       Characterizing Mal'cev conditions, \AU\ 3 (1973),
           351--397.
\bibitem{wtaylor-fsv}
 \Dash
      The fine spectrum of a variety, \AU\ 5 (1975), 263--303.
\bibitem{wtaylor-vohl}
 \Dash
      Varieties obeying homotopy laws, \CANAD\ 29 (1977),
      498--527.
\bibitem{wtaylor-cots}
 \Dash
 {\em The clone of a topological space}, Volume 13 of
     {\em Research and Exposition in Mathematics}, 95 pages.
     Heldermann Verlag, 1986.
\bibitem{wtaylor-gcg}
  \Dash
      {\em The Geometry of Computer Graphics}, Wadsworth and
      Brooks-Cole, Pacific Grove, CA, 1992.
\bibitem{wtaylor-sae}
  \Dash
       Spaces and equations, \FM\  164 (2000), 193--240.
        %
\bibitem{wtaylor-eri}
  \Dash
       Equations on real intervals, \AU\/\  55 (2006), 409--456.
\bibitem{wtaylor-asi}
   \Dash
         Approximate satisfaction of identities, 98 pp., 2010. See 
         \url{http://arxiv.org/abs/1504.01165}
\bibitem{wtaylor-class} 
       \Dash 
         Classification of finite-dimensional compact topological
         algebras. See \url{http://arxiv.org/abs/1402.3734}
\bibitem{vanmill}
J. van Mill,
          A topological group having no homeomorphisms
          other than translations. \TAMS\ 280 (1983),  491--498.
\bibitem{wallace}
A. D. Wallace,
        The structure of topological semigroups, \BAMS\ 61
        (1955), 95--112.
\end{thebibliography}
